\documentclass[12pt]{amsart}
\usepackage{amsmath}
\usepackage{amssymb}
\usepackage{graphicx}
\usepackage{graphics}
\usepackage{amsthm}
\usepackage{amscd}
\usepackage{color}
\usepackage{amsfonts}
\usepackage{fancyhdr}
\newtheorem{theorem}{Theorem}[section]
\newtheorem{mainthm}{Theorem}
\newtheorem*{theorem*}{Theorem}
\newtheorem{corollary}[theorem]{Corollary}
\newtheorem{proposition}[theorem]{Proposition}
\newtheorem{lemma}[theorem]{Lemma}
\newtheorem{Remark}[theorem]{Remark}
\newtheorem*{definition*}{Definition}

\newtheorem{claim}[theorem]{Claim}
\newtheorem{definition}[theorem]{Definition}
\newtheorem*{Question*}{Question}
\newtheorem{Question}{Question}
\newtheorem{Conjecture}{Conjecture}

\def\Z{\mathbb{Z}}

\def\R{\mathbb{R}}

\def\T{\mathbb{T}}

\def\norm #1{\Vert \,#1\, \Vert\,}
\makeatletter

\newcommand{\Rmnum}[1]{\expandafter\@slowromancap\romannumeral #1@}
\makeatother

\def\ud{\mathrm{d}}
\def\diff {\operatorname{Diff}}
\def\dim{\operatorname{dim}}
\def\homeo{\operatorname{Homeo}}
\def\Orb{\operatorname{Orb}}
\def\um{\operatorname{m}}
\def\Id{\operatorname{Id}}

\def\La{\Lambda}

\def\cA{\mathcal{A}}    
  \def\cH{\mathcal{H}}  \def\cT{\mathcal{T}}
    \def\cU{\mathcal{U}}
   \def\cP{\mathcal{P}} \def\cV{\mathcal{V}}
    
\def\cF{\mathcal{F}}  \def\cL{\mathcal{L}} \def\cR{\mathcal{R}}

\keywords{Partial hyperbolicity, dynamical coherence, complete foliation, accessible boundary, Anosov flow}


\begin{document}
\vspace{-2cm}
\title[Partially hyperbolic diffeomorphisms with neutral center ]
{Partially hyperbolic diffeomorphisms with one-dimensional neutral center on
3-manifolds }
\author{Jinhua Zhang }
\vspace{-2cm}
\maketitle
\begin{abstract}We prove that for any partially hyperbolic diffeomorphism   with one dimensional neutral center on a 3-manifold, the center stable and center unstable foliations are complete; moreover, each leaf of center stable and center unstable foliations  is a cylinder, a M$\ddot{o}$bius band or a plane.

 Further properties of the Bonatti-Parwani-Potrie type of partially hyperbolic diffeomorphisms  are studied. Such examples   are obtained by composing the time $m$-map (for $m>0$ large) of a non-transitive Anosov flow $\phi_t$ on an orientable 3-manifold with Dehn twists along some transverse tori, and the examples are partially hyperbolic with one-dimensional neutral center. We prove that the center foliation gives a topologically Anosov flow which is topologically equivalent to $\phi_t$. We also prove that for the precise example constructed by Bonatti-Parwani-Potrie, the center stable and center unstable foliations are robustly complete.

\end{abstract}


\section{Introduction}
\subsection{Our setting}
In 1970s, M. Brin and Y. Pesin ~\cite{BP} proposed a  notion called \emph{partial hyperbolicity}. A diffeomorphism is called partially hyperbolic if the tangent bundle of the manifold splits into three invariant bundles: one of which is uniformly contracting under the dynamics, another is uniformly expanding, and the center is intermediate.

One of the main topics on partially hyperbolic systems is the classification according to different properties.  This field is the intersection of topology and dynamical systems,  and many important projects are proposed by  F. Rodriguez Hertz, J. Rodriguez Hertz and R. Ures in a series of papers and talks.
Extending a conjecture by E. Pujals (formalized in~\cite{BW}),  F. Rodriguez Hertz, J. Rodriguez Hertz and R. Ures ~\cite{HHU6} proposed the following:
\begin{Conjecture}\label{conjecture.hhu}
Any dynamically coherent partially hyperbolic   diffeomorphism on a 3-manifold is, up to finite iterations and finite lifts,  leaf conjugated to one of the following three models:
\begin{itemize}
\item linear Anosov diffeomorphism on $\T^3$;
\item  time one map of an Anosov flow;
\item skew products over  linear Anosov diffeomorphisms on torus.
\end{itemize}
\end{Conjecture}
We remark that it is necessary to consider finite  iterations and  lifts, see the examples in ~\cite[Section 4]{BW}.

For this conjecture, some partial  results are obtained, see for instance ~\cite{Bo},~\cite{BW},
 \newline
 ~\cite{Ca}, ~\cite{HaPo0}, \cite{HaPo}, \cite{G}.  In~\cite{BW}, a notation called ~\emph{completeness} was proposed for the structure of invariant foliations of a \emph{dynamically coherent} partially hyperbolic diffeomorphism and the completeness of  invariant foliations played an important role in attacking this conjecture. We remark that for all the three models in Conjecture~\ref{conjecture.hhu}, their invariant foliations are complete.

Recently, C. Bonatti, K. Parwani and R. Potrie ~\cite{BPP}  gave a mechanism to build new partially hyperbolic diffeomorphisms which are counter examples to the conjecture above.   The new partially hyperbolic diffeomorphism is obtained by composing the time $n$-map of a specific non-transitive Anosov flow with a Dehn twist along a transverse torus.
 Then in ~\cite{BZ}, it is shown that such construction can be made to any non-transitive Anosov flow.
  One can ask: to what extent  the properties of the three models  in the conjecture are preserved by the new partially hyperbolic diffeomorphisms? To be precise:
\begin{Question*} For the new partially hyperbolic diffeomorphisms in ~\cite{BPP,BZ}, is every center stable  (resp. center unstable) leaf either a cylinder or a plane?
Are  the center stable  and center unstable foliations  complete?  Furthermore, what is the relation between the new partially hyperbolic diffeomorphism and the Anosov flow used for building it?
\end{Question*}
In this paper, we give answers to the questions above.

\subsection{Statement of the results}
Let $M$ be a closed 3-manifold. We denote by $\cP\cH(M)$ the set of partially hyperbolic diffeomorphisms with the splitting of the form
$$TM=E^s\oplus E^c\oplus E^u, \textrm{ where $\dim(E^s)=\dim(E^c)=1$.}$$
Given $f\in\cP\cH(M)$, one says that  $f$ has \emph{ neutral behavior along  $E^c$ or neutral center} if there exists a constant $K>1$ such that
$$\frac 1K\leq\norm{Df^n|_{E^c(x)}}\leq K, \textrm{ for any $x\in M$ and any $n\in\Z$.}$$
 By Theorem 7.5 of ~\cite{HHU2}, one has that $f$ is \emph{dynamically coherent}, that is, there exist $f$-invariant foliations $\cF^{cs}$ and $\cF^{cu}$ tangent
 to $E^s\oplus E^c$ and $E^c\oplus E^u$ respectively. Moreover, one has that the invariant foliations $\cF^{cs}$ and $\cF^{cu}$ are \emph{plaque expansive}.
 It is well known that $E^s$ and $E^u$ are uniquely integrated into two invariant foliations $\cF^{ss}$ and $\cF^{uu}$. By Remark 3.7  of ~\cite{HHU3}, the neutral behavior along the center implies that the center distribution $E^c$ is also  uniquely integrable.

  Given a partially hyperbolic diffeomorphism $f$, assume that $f$ is dynamically coherent,  then one can get an $f$-invariant foliation $\cF^c$ tangent to $E^c$. We denote  $\cF^{ss}(\cF^c(x)):=\cup_{y\in\cF^c(x)}\cF^{ss}(y)$. We say that \emph{ the center stable foliation is complete}  if for any $x\in M$, we have that
  $$\cF^{ss}(\cF^c(x))=\cF^{cs}(x).$$
  Although we don't have the examples of dynamically coherent partially hyperbolic diffeomorphisms whose invariant foliations are not  complete, we cannot rule out this possibility.

  Our first result is the following:
\begin{mainthm}\label{thm.plane and cylinder} Let $M$ be a closed 3-manifold and $f\in\cP\cH(M)$.
Assume that $f$ has neutral behavior along  the center, then we have the followings:
\begin{itemize}
\item[--] the center stable and center unstable foliations are complete;
\item[--] every center stable (resp. center unstable) leaf  is   a plane,  a M$\ddot{o}$bius band or a cylinder.
Moreover, a center stable (resp. center unstable) leaf is a cylinder or a M$\ddot{o}$bius band  if and only if this leaf contains a compact center leaf.
\end{itemize}
\end{mainthm}
\begin{Remark}\label{r.robust of second}
  By ~\cite[Theorem 7.5]{HHU2}, one has that both the  center stable and center unstable foliations are plaque expansive.
Hence, there exists a  small neighborhood $\cU$ of $f$ such that each  $g\in\cU$ is dynamically coherent and
 every center stable (center unstable) leaf is   a plane,  a M$\ddot{o}$bius band or a cylinder.
\end{Remark}
By Remark ~\ref{r.robust of second}, the second item of Theorem~\ref{thm.plane and cylinder} is a robust property. However, we don't know if the first property is  robust.
\begin{Question} Does there exist a small neighborhood $\cV$ of $f$ such that for any $g\in\cV$, the center stable and center unstable foliations of $g$ are complete?
\end{Question}

Let $\phi_t$ be a smooth non-transitive  Anosov flow on an orientable 3-manifold $M$. Consider a smooth \emph{Lyapunov function}
 $\cL:M\mapsto \R$ of the flow $\phi_t$ (for definition see Section~\ref{s.anosov flow}). Let $\{c_1,\cdots,c_m\}$ be the values of $\cL$ on the hyperbolic basic sets. We say that $\cL^{-1}(c)$ is a \emph{wandering regular level of  $\cL$} if $c$ is a regular value of $\cL$ and $c$ is in $\cL(M)\backslash\{c_1,\cdots,c_m\}$.
 Then the  wandering regular level $\cL^{-1}(c)$  consists of  finite pairwise disjoint tori transverse
  to the Anosov flow, see for instance \cite{Br}.

  Now, we define the set of partially hyperbolic diffeomorphisms that we consider. Given an orientable 3-manifold $M$ and let $\phi_t$ be a non-transitive Anosov flow on $M$, we denote by
$\mathcal{PH}_{\phi_t}(M)\subset \mathcal{PH}(M)$ the set of diffeomorphisms such that for each $f\in\mathcal{PH}_{\phi_t}(M)$, one has that
\begin{itemize}
\item $f$ is partially hyperbolic with one dimensional center;
\item there exist $\tau>0$ and a family of tori $\{T_1,\cdots,T_k\}$ contained in a wandering regular level of a smooth Lyapunov function of $\phi_t$  such that
$$f=\psi_1\circ\cdots\circ\psi_k\circ\phi_{\tau},$$
where  $\psi_i$ is a Dehn twist along $T_i$ and is supported  in  $\{\phi_t(T_i)\}_{t\in(0,\tau)}$.
\end{itemize}

\begin{Remark}
It is shown in ~\cite{BZ} that for each smooth non-transitive Anosov flow $\phi_t$ on $M$, one always has that  $\mathcal{PH}_{\phi_t}(M)$ is  non-empty,  see \cite[Theorem 7.1 and Lemma 7.4]{BZ}.
\end{Remark}
 \begin{mainthm}~\label{thm.conjugate to Anosov flow}
 Let $\phi_t$ be a non-transitive Anosov flow on an orientable 3-manifold $M$. For any  $f\in\mathcal{PH}_{\phi_t}(M)$, one has that
 \begin{itemize}
 \item[--] the diffeomorphism $f$ has neutral center;

\item[--] We denote by $\cF^c$ the center foliation of $f$, then there exist a continuous flow $\{\theta_t\}_{t\in \mathbb{R}}:M\mapsto M$ and a homeomorphism $h:M\mapsto M$ such that for any $x\in M$, one has

 \begin{itemize}

 \item[--] $$\Orb(x, \theta_t)=\cF^c(x);$$
 \item[--] $$h(\Orb(x,\theta_t))=\Orb(h(x),\phi_t)\textrm{ and }h(\Orb^+(x,\theta_t))=\Orb^+(h(x),\phi_t);$$
 \end{itemize}

 \end{itemize}

 \end{mainthm}

 Now, we discuss the particular example $f_b$ built in ~\cite{BPP}. In order to state the further  properties of $f_b$ that we get, we need to recall some terminology to give the statement of our result.
In ~\cite[Section 4]{BPP}, the authors firstly  build a smooth non-transitive Anosov flow $\psi_t$ on a 3-manifold $N$ such that
 \begin{itemize}
 \item the non-wandering set consists of  one attractor and one repeller.
 \item there exist two transverse tori  and each  orbit has no return on the union of these two tori.
 \item the foliations induced by the stable and unstable foliations of the Anosov flow on each transverse torus consist of two Reeb components;
 \item the union of these two tori separates the manifold into two connected components which contain the attractor and the repeller respectively.
 \end{itemize}
 The partially hyperbolic diffeomorphism $f_b$ in ~\cite{BPP} is obtained by composing a Dehn twist along one transverse torus with $\psi_n$ (for $n>0$ large), and $f_b$ has one dimensional neutral center (see \cite[Lemma 9.1]{BPP}).

  By Proposition 1.9 in ~\cite{BZ}, the two transverse tori for $\psi_t$
  are contained in a wandering regular level of a smooth Lyapunov function of $\psi_t$.
    Hence Theorem~\ref{thm.conjugate to Anosov flow} can be applied to $f_b$. This shows that  the action of $f_b$ on the space of center leaves   just ``permutes" the center leaves of the Anosov flow on each center stable (center unstable) leaf without changing  the structure of
 the invariant foliations of the Anosov flow $\psi_t$.
 Combining with \cite[Theorem 9.6]{BPP}, one has the following corollary:
    \begin{corollary} There exist a 3-manifold $M$ and a dynamically coherent partially hyperbolic diffeomorphism $f$ on $M$ such that
    \begin{itemize}
    \item the manifold $M$ supports Anosov flows;
    \item the center foliation of $f$ defines a continuous flow which is topologically equivalent to an Anosov flow;
\item  there is no lifts or iterations of $f$ that is leaf conjugate to the time one map of an Anosov flow.
    \end{itemize}
    \end{corollary}

  By Theorem 7.1 and Lemma 7.4 in ~\cite{BZ}, one can apply Theorems ~\ref{thm.plane and cylinder} and ~\ref{thm.conjugate to Anosov flow} to the new examples in ~\cite{BZ}, and one gets that  their center stable (resp. center unstable) leaves are  planes,  cylinders, or M$\ddot{o}$bius bands;
moreover, their center foliations are topologically Anosov. One can ask a more ambitious question:
\begin{Question}
Given a dynamically coherent partially hyperbolic diffeomorphism $f\in\cP\cH(M)$.  Up to finite lifts and finite iterations, does one of the followings hold:
\begin{itemize}
\item  $f$ is leaf conjugate to a linear Anosov diffeomorphism on $\mathbb{T}^3$;
\item or $f$ is leaf conjugate to a skew product over a linear Anosov diffeomorphism on $\mathbb{T}^2$;
\item or the center foliation of $f$ defines a continuous flow which is topologically equivalent to an Anosov flow?
\end{itemize}
\end{Question}

 By Theorem~\ref{thm.plane and cylinder}, one also has that for the examples in ~\cite{BZ}, the center stable and center unstable foliations are complete. Then one can ask:
  \begin{Question} For the partially hyperbolic diffeomorphisms with neutral center built in ~\cite{BZ}, are the center stable and center unstable foliations  robustly complete?
 \end{Question}
 We don't have the answer for the question above, but one can get the robust completeness of the center stable and center unstable foliations for the particular  example $f_b$ in ~\cite{BPP}.

  \begin{proposition}\label{thm.conjugation} There exists a $C^1$ neighborhood $\cU$ of $f_b$ such that for any $g\in\cU$, one has that the center stable and  center unstable foliations of $g$ are   complete.
 \end{proposition}

\begin{Remark}
The proof of Proposition ~\ref{thm.conjugation} depends on the fact that the non-wandering set of $\phi_t$ consists of one attractor and one repeller.
\end{Remark}

\section{Preliminary}
In this section, we collect the notations and results that we need.

Let $f$ be a diffeomorphism on a compact manifold $M$,  and recall that $f$ is \emph{partially hyperbolic}, if there exist a $Df$ invariant splitting $TM=E^s\oplus E^c\oplus E^u$ and a positive integer $N$ such that for any $x\in M$, one has
$$\norm{Df^N|_{E^s(x)}}<\min\{\um(Df^N|_{E^c(x)}),1\}\leq\max\{\norm{Df^N|_{E^c(x)}},1\}
<\um(Df^N|_{E^u(x)}).$$
\subsection{ Dynamical coherence, plaque expansiveness and completeness}
\begin{definition}
Let $f$ be a partially hyperbolic diffeomorphism. We say that $f$ is \emph{cs (resp. cu)-dynamically coherent}, if there exists an $f$-invariant foliation $\cF^{cs}$ (resp. $\cF^{cu}$)  tangent to $E^{s}\oplus E^c$ (resp. $E^{c}\oplus E^u$). In particular, $f$ is \emph{dynamically coherent}, if $f$ is both cs-dynamically coherent and cu-dynamically coherent.
\end{definition}
A partially hyperbolic diffeomorphism might  not be  dynamical coherent even if the center dimension is one (see for instance~\cite{HHU4}).

Let $\cF$ be an $f$-invariant foliation. We denote by $\cF(x)$ the $\cF$-leaf through the point $x$ and by $\cF_\epsilon(x)$ the $\epsilon$-neighborhood of $x$ in the leaf $\cF(x)$.
A sequence of points $\{x_n\}_{n\in\Z}$ is called \emph{an $\epsilon$ pseudo orbit with respect to $\cF$}, if $f(x_n)$ belongs to $\cF_{\epsilon}(x_{n+1})$ for any $n\in\Z$.
\begin{definition} Given an $f$-invariant foliation $\cF$.
We say that  $\cF$ is \emph{plaque expansive}, if there exists $\epsilon>0$ satisfying the following:
if $\{x_n\}_{n\in\Z}$ and $\{y_n\}_{n\in\Z}$ are two $\epsilon$-pseudo orbits with respect to $\cF$ and if one has $\ud(x_n,y_n)<\epsilon$ for any $n\in\Z$,
then $x_n$ and $y_n$ belong to a common $\cF$-leaf for any $n\in\Z$.
\end{definition}

By Theorem 7.5 and Corollary 7.6  in ~\cite{HHU2},  a partially hyperbolic diffeomorphism with neutral center is dynamically coherent; moreover, the center, center stable and center unstable foliations are plaque expansive. By Theorem 7.1 in ~\cite{HPS} one has that the plaque expansiveness in this setting is a robust property and implies the structure stability of the invariant foliation (ie. leaf conjugacy). In ~\cite{PS}, the authors prove that if the center foliation is plaque expansive, then the leaf conjugacy for the center foliation is also the leaf conjugacy for  the center stable and center unstable foliations. To summarize, one has the following result:
\begin{theorem}\label{r.preserves all the leaves}
Let $f$ be a partially hyperbolic diffeomorphism. If $f$ has one dimensional neutral center, there exists a $C^1$ small neighborhood $\cU$ of $f$ such that for any $g\in\cU$, one has the following properties:
 \begin{itemize}
\item $($dynamical coherence$)$  $g$ is dynamically coherent;
\item $($plaque expansive$)$  the center, center stable and center unstable foliations $g$ are plaque expansive;
\item $($leaf conjugacy$)$ there exists a homeomorphism $h_g:M\mapsto M$ such that for any point $x\in M$ and  $i=c,cs,cu$, one has that
$$h_g(\cF^{i}_g(x))=\cF^{i}_f(h_g(x))\textrm{ and }h_g(g(\cF_g^{i}(x)))=f(h_g(\cF_g^{i}(x))).$$
 \end{itemize}
\end{theorem}
\begin{Remark} The homeomorphism $h_g$ tends to identity as $g$ tends to $f$.
\end{Remark}

Let $f$ be a dynamically coherent partially hyperbolic diffeomorphism. A point $y$ is called \emph{an accessible boundary point} with respect to  $\cF^{ss}(\cF^c(x))$ for some $x\in M$, if there exists a $C^1$ curve $\sigma:[-1,0]\mapsto M$ tangent to the center bundle such that
 $$\sigma([-1,0))\subset \cF^{ss}(\cF^c(x)) \textrm{ and }  \sigma(0)=y\notin  \cF^{ss}(\cF^c(x)) .$$
 The set of accessible boundary points  with respect to    $\cF^{ss}(\cF^c(x))$ is called \emph{accessible boundary} with respect to     $\cF^{ss}(\cF^c(x))$.

 With the notations above,  one has the following result due to ~\cite{BW}:
\begin{proposition}\label{p.boundary leaf} The accessible boundary with respect to  $\cF^{ss}(\cF^{c}(x))$ is saturated by strong stable leaves.
\end{proposition}
We will call each strong stable leaf in the accessible boundary  with respect to  $\cF^{ss}(\cF^{c}(x))$ as \emph{a boundary leaf with respect to  $\cF^{ss}(\cF^{c}(x))$ or a boundary leaf }for simplicity.
\subsection{Existence of compact leaves}

In 1965, S. Novikov gave a criterion for the existence of  compact leaves of codimension one foliations on closed 3-manifolds.
\begin{theorem}\label{thm.Novikov}\cite{N}
Let $\cF$ be a codimension one foliation  on a  3-manifold $M$. $\cF$ has a compact leaf,  if one of the followings is satisfied:
\begin{itemize}
\item there exists a null-homotopy closed  transversal  for $\cF$;
\item there exists a non-null homotopic closed path in a $\cF$-leaf which is null homotopy in $M$.
\end{itemize}
\end{theorem}
With the help of Novikov's theorem, \cite{HHU5} proves the non-existence of compact leaf for center stable (center unstable) foliation. More precisely,
\begin{theorem}\cite[Theorem 1.1]{HHU5}\label{thm.no compact leaf} Let $f$ be a partially hyperbolic diffeomorphism on a 3-manifold $M$ with the splitting $TM=E^s\oplus E^c\oplus E^u$.
Assume that $f$ is cs-dynamically coherent, 
 then the center stable foliation $\cF^{cs}$ has no compact leaves.
\end{theorem}
In our context, instead of Theorem~\ref{thm.no compact leaf},
we can also  use the following result\footnote{I thank Rafael Potrie for pointing out this fact.}:
\begin{theorem}\cite{HHU1} Assume that $f$ is partially hyperbolic diffeomorphism on a 3-manifold $M$ and admits a torus tangent to $E^s\oplus E^c$,
then $M$ fibers over $S^1$ with torus fiber.
\end{theorem}
In fact, the dynamical coherent partially hyperbolic diffeomorphisms on 3-manifolds fibering $S^1$ over torus do not admit compact center stable or center unstable leaves,   see ~\cite{HaPo0, HaPo}.
\subsection{Anosov flow}\label{s.anosov flow}

A non-singular $C^1$ vector field $X$ on a closed manifold $M$ is called \emph{an Anosov vector field}, if there exists a splitting
  $TM=E^s\oplus<X>\oplus E^u$ which is invariant under the tangent flow $\Phi^X_t$ of  flow $X_t$ such that $E^s$ is uniformly contracting and $E^u$ is uniformly expanding under $\Phi^X_t$. The flow generated by an Anosov vector field is called as \emph{Anosov flow}.

It is well known that the Anosov flow is structurally stable, that is, the flows generated by the  vector fields in a $C^1$ small neighborhood of  an Anosov vector field are topologically equivalent to each other. The topological  equivalence of two flows is defined as follow:
\begin{definition}
Let $\varphi_t: M\mapsto M$ and $\theta_t:N\mapsto N$ be  two continuous flows. We say that \emph{$\varphi_t$ is topologically equivalent  to $\theta_t$}, if there exists a homeomorphism $h:M\mapsto N$ preserving the orientation of the flows and sending the orbits of  the flow $\varphi_t$ to the orbits of the flow $\theta_t$, that is, for any $x\in M$, one has
$$h(\Orb(x,\varphi_t))=\Orb(h(x),\theta_t) \textrm{ and  } h(\Orb^+(x,\varphi_t))=\Orb^+(h(x),\theta_t).$$
\end{definition}

Different from the Anosov diffeomorphisms on surfaces, the Anosov flows on 3-manifolds might  be non-transitive, see for instance~\cite{FW}.
Given an Anosov flow $\phi_t$ on $M$,  a  smooth \emph{Lyapunov function} for $\phi_t$ is a smooth function $\cL: M\mapsto\R$ such that
 \begin{itemize}
 \item $\cL$ is not increasing along every orbit;
 \item  $\cL$ is strictly decreasing along an orbit  if and only if  this orbit is in the wandering domain.
 \end{itemize}
 For a smooth Anosov flow, there always exist smooth Lyapunov functions, see \cite[Page 18]{S}.

 Recall that an embedded torus defined by $e: \T^2\mapsto M$ is \emph{incompressible}, if the  induced map  $e_{*}:\pi_1(T)\mapsto\pi_1(M)$ is injective.  With the help of Lyapunov functions,
 for non-transitive Anosov flows, one can separate  the hyperbolic basic sets by finite incompressible  transverse tori, see ~\cite{Br}.
\subsection{Dehn twists}
Now, we give the definition of Dehn twists on 3-manifolds.

\begin{definition}Let $T$ be an embedded torus on a 3-manifold $M$. We say that a diffeomorphism $\Psi:M\mapsto M$ is a \emph{Dehn twist along the torus $T$ },
if there exists an orientation preserving  diffeomorphism $\varphi:\T^2\times[0,1]\hookrightarrow M$ such that
 \begin{itemize}
 \item $\varphi(\T^2\times\{0\})=T$,
 \item $\Psi$ is identity in the complement of $\varphi(\T^2\times[0,1])$;
 \item under the coordinate of $\varphi$, one has that
  \begin{enumerate}
  \item $\varphi^{-1}\circ\Psi\circ\varphi:\T^2\times[0,1]\mapsto\T^2\times[0,1]$ is of the form
$$\varphi^{-1}\circ\Psi\circ\varphi(x,t)=(\phi_t(x),t),\textrm{ for any $(x,t)\in\T^2\times[0,1]$;}$$
\item $\phi_t$ equals to identity when $t$ is close to 0 or 1;
\item for each $x\in\T^2$, the closed path $\{\phi_t(x)\}_{t\in[0,1]}$ is non-null homotopy.
\end{enumerate}
\end{itemize}
\end{definition}

\section{The topological structure of the center stable and center unstable
foliations: Proof of Theorem ~\ref{thm.plane and cylinder}}

Let $f$ be a partially hyperbolic diffeomorphism on a 3-manifold $M$,  exhibiting neutral behavior along one dimensional center.  Since it has already been proven that $f$ is dynamically coherent, we denote   the center stable and center unstable foliations as $\cF^{cs}$ and $\cF^{cu}$ respectively.
We will first show that $\cF^{cs}$ and $\cF^{cu}$ are complete, then we give the description of their  leaves.
\proof [Proof of Theorem~\ref{thm.plane and cylinder}]

By the uniform transversality of strong stable direction and center direction restricted to every center stable leaf, there exists $\delta>0$ such that for any point $x\in M$,  the $\delta$ neighborhood of the leaf $\cF^{c}(x)$ in the center stable leaf $\cF^{cs}(x)$  is contained in $\cF^{ss}(\cF^{c}(x))$.

We will prove the completeness of $\cF^{cs}$ by contradiction. Assume that $\cF^{cs}$ is not complete, then by Proposition~\ref{p.boundary leaf},
there exists a point $x\in M$ such that $\cF^{ss}(\cF^{c}(x))$ has   a boundary leaf $\cF^{ss}(y)$ for some $y\in M$ (there might be infinitely many boundary leaves).  By the invariant property of the center  and strong stable foliations, we have that $\cF^{ss}(f^{n}(y))$ is a boundary leaf with respect to  $\cF^{ss}(\cF^{c}(f^{n}(x)))$, for any integer $n\in\Z$.  By the choice of $\delta$, one has  that when it is restricted to  the center stable leaf $\cF^{cs}(f^{n}(x))$, the strong stable leaf $\cF^{ss}(f^{n}(y))$ is $\delta$ away from the center leaf $\cF^{c}(f^n(x))$.

  By the definition of  boundary leaves, there exists a $C^1$-curve $\sigma: [0,1]\mapsto M$ such that
  $$\sigma(t)\subset \cF^c(y),  \sigma(0)=y \textrm{  and  }\sigma((0,1])\subset\cF^{ss}(\cF^c(x)).$$
  Up to  shrinking $\sigma$, we can assume that the length $\ell(\sigma)$ of $\sigma$ is strictly less than $\frac{\delta}{4K^2}$, where $K>1$ is the number satisfying
  $$\frac{1}{K}\leq \norm{Df^n|_{E^c(p)}}\leq K, \textrm{ for any $n\in\Z$ and any point $p\in M$}.$$
   Since $\sigma(1)$ is on the strong stable manifold of a point $z\in\cF^c(x)$, there exists an integer $m$ large enough such that $f^{m}(\sigma(1))$ is in the $\frac \delta2$ neighborhood of $f^{m}(z)$ with respect to the distance on  the center stable leaf $\cF^{cs}(f^{m}(x))$. Since $f^{m}(y)$ is still on the boundary leaf $\cF^{ss}(f^{m}(y))$, we have that the length
  $$\ell(f^m(\sigma))>\frac\delta2.$$
  On the other hand, we have the estimate
  $$\ell(f^m(\sigma))\leq \max_{p\in M}\norm{Df^{m}|_{E^c(p)}}\cdot\ell(\sigma)<\frac{\delta}{4K}<\frac{\delta}{2},$$
  a contradiction. This proves the first item of Theorem~\ref{thm.plane and cylinder}.

Let $\tilde{M}$ be the universal cover of $M$. The metric on $\tilde{M}$ is the pull back of the metric on $M$ by the covering map. We denote by $\tilde{\cF^{i}}$   the lift of $\cF^{i}$ on the universal cover $\tilde{M}$, for $i=ss,cs,c,cu,uu$.

 To prove the second item, we need the following lemma:
 \begin{lemma}\label{l.lifted to plane} For any $x\in\tilde{M}$, the lifted leaves $\tilde{\cF}^{cs} (x)$ and $\tilde{\cF}^{cu}(x)$ are planes.
 \end{lemma}
 \proof For any $x\in\tilde{M}$,  the lifted leaf $\tilde{\cF}^{cs} (x)$ is a two dimensional manifold without boundary; to prove that it is a plane, we only need to show that its fundamental group is trivial.  Assume  that there exists a closed curve $\gamma$ in $\tilde{\cF}^{cs}(x)$ which is non-null homotopy in the leaf. Since  $\gamma$ is null homotopy in $\tilde{M}$,  then  the projection of $\gamma$ on $M$ is null homotopy in $M$ and is non-null homotopy in a $\cF^{cs}$-leaf.  By Theorem~\ref{thm.Novikov}, the foliation $\cF^{cs}$ has a compact leaf,  which   contradicts to   Theorem~\ref{thm.no compact leaf}.

 Analogously, one can show that $\tilde{\cF}^{cu} (x)$ is also a plane.
 \endproof
 \begin{claim}\label{c.trivial foliated} The lifted foliations  $\tilde{\cF}^{cs }$ and $\tilde{\cF}^{cu} $ are complete, that is, they are trivially bi-foliated by $\tilde{\cF}^{ss} $ and $\tilde{\cF}^{c} $.
 \end{claim}
 \proof
  Let $\tilde{f}$ be a lift of $f$, then $\tilde{f}$ is a partially hyperbolic diffeomorphism with one dimensional  neutral  center  and whose invariant foliations are the lifts of the invariant foliations of $f$. Hence, the argument for $f$ applies for $\tilde{f}$.
 \endproof

 For any compact center leaf $\gamma$ of $f$, we prove the following:
  \begin{lemma} The center stable leaf contains $\gamma$ is either a cylinder or a M$\ddot{o}$bius band.
  \end{lemma}
  \proof   Up to taking a double cover of the manifold, we can assume that the strong stable bundle is orientable and we give it an orientation.
By Theorem~\ref{thm.no compact leaf}, the leaf $\cF^{cs}(\gamma)$ is not compact.

  If every strong stable leaf through a point on $\gamma$ intersects $\gamma$ only once, then one claims that   the center stable leaf of $\gamma$ is a cylinder. On the universal cover,  a lift $\tilde{\cF}^{cs}(\gamma)$ of $\cF^{cs}(\gamma)$ is a plane. By  completeness and the fact that each strong stable leaf through $\gamma$ intersects $\gamma$ only once, one has that the lift of $\gamma$ in $\tilde{\cF}^{cs}(\gamma)$ has only one connected component. Let $\Gamma\subset\pi_1(M)$ be the subgroup which keeps $\tilde{\cF}^{cs}(\gamma)$ invariant, then $\cF^{cs}(\gamma)$ is the quotient of $\tilde{\cF}^{cs}(\gamma)$ by the action of $\Gamma$ and $\Gamma$ is the fundamental group of $\cF^{cs}(\gamma)$. Let $\tilde{\gamma}$ be the lift of $\gamma$ in $\tilde{\cF}^{cs}(\gamma)$, then $\tilde{\gamma}$ is homeomorphic to $\mathbb{R}$ and each element of $\Gamma$ keeps $\tilde{\gamma}$ invariant. Since each non-trivial element of $\Gamma$ has no fixed points, each non-trivial element acting  on $\tilde{\gamma}$ has no fixed points. To summarize, one has that $\Gamma$ is isomorphic to a subgroup of $\homeo(\mathbb{R})$, whose non-trivial element has no fixed points.   H$\ddot{o}$lder theorem~\cite{Ho} asserts that any group acting freely  on $\R$ is Abelian. Hence $\Gamma$ is Abelian, then $\cF^{cs}(\gamma)$ can only be a cylinder.

  If not, there exists a strong stable leaf intersects $\gamma$ at least twice. Consider the universal cover $\tilde{M}$ and the lift $\cP$ of the leaf $\cF^{cs}(\gamma)$ which is a plane, then  there exist two center leaves which are  the  lifts of $\gamma$ on $\cP$; by the completeness of the center stable foliation on the universal cover, one has that for every point on $\gamma$, the strong stable curve through this point positively goes back to $\gamma$ after some uniform finite length, which implies that the center stable leaf $\cF^{cs}(\gamma)$ is a closed surface, a contradiction.
   \endproof
 The following lemma ends the proof of Theorem~\ref{thm.plane and cylinder}:
 \begin{lemma}\label{l.plane 3.3} Any center stable leaf which  contains no  compact center leaves is  a plane.
 \end{lemma}
\proof Let $\cF^{cs}(x)$ be  a center stable leaf which   contains no  compact center leaves. We will first use the argument from ~\cite{BW} to prove that $\cF^{cs}(x)$ is either a plane or a cylinder,  then we show that $\cF^{cs}(x)$ can only be  a plane.

Consider a lift $\tilde{\cF}^{cs}(y)$ of  $\cF^{cs}(x)$  and  let $\Gamma$ be the subgroup of $\pi_{1}(M)$ which keeps the leaf $\tilde{\cF}^{cs}(y)$ invariant, then  $\cF^{cs}(x)$ is the quotient of $\tilde{\cF}^{cs}(y)$ by the action of $\Gamma$.
By  Claim~\ref{c.trivial foliated}, we have that the space of center leaves   in $\tilde{\cF}^{cs}(y)$ is a   real  line, as well as the space of strong stable leaves in $\tilde{\cF}^{cs}(y)$.  Hence, $\Gamma$ induces two actions on these two spaces respectively  and the action of $\Gamma$ is a sub-action of Cartesian product of these two actions.

The action of $\Gamma$  on the space of center leaves corresponds to  a subgroup of $\homeo(\R)$.  Moreover every non-trivial element of  $\Gamma$ acts on the space of center leaves without fixed points and  preserving the orientation,  otherwise there exists a non-trivial element of  $\Gamma$ keeping a center leaf invariant which implies that $\cF^{cs}(x)$ has a compact center leaf. Similar argument applies to the action of $\Gamma$ on the space of the strong stable leaves,  proving that this action is orientation preserving and has no fixed points.    Recall that H$\ddot{o}$lder theorem ~\cite{Ho} asserts that any group acting freely on $\R$ is Abelian, hence these two actions are Abelian actions.  As a consequence, the action of  $\Gamma$ is Abelian and is orientation preserving, which  implies that $\cF^{cs}(x)$ is an orientable  two dimensional  manifold whose fundamental group is Abelian. Once again, by Theorem~\ref{thm.no compact leaf}, we have that $\cF^{cs}(x)$ is either a cylinder or a plane.

By assumption, one has that $\cF^{cs}(f^n(x))$ contains no compact center leaves for any integer $n\in\Z$.
We will prove,  by contradiction, that $\cF^{cs}(x)$ is a plane. Assume that  $\cF^{cs}(x)$ is a cylinder, then one has the following result:
\begin{claim} The center leaf $\cF^c(x)$ intersects $\cF^{ss}(x)$ at least twice.
 \end{claim}
 \proof Let $y$ be a lift of $x$. Since  leaf $\tilde{\cF}^{cs}(y)$ is complete and the group $\Gamma$ is non-trivial, given a non-trivial element $\varphi\in\Gamma$,
 then $\varphi(\tilde{\cF^c}(y))$ is a leaf different from $\tilde{\cF^c}(y)$ (otherwise, $\cF^c(x)$ is compact). As a consequence,  $\varphi(\tilde{\cF^c}(y))$ intersects $\tilde{\cF}^{ss}(y)$ in a point different from  $y$ and $\varphi(y)$, hence $\cF^c(x)$ intersects $\cF^{ss}(x)$ at least twice.
 \endproof
One can  take a point $z\in \cF^c(x)\cap\cF^{ss}(x)\backslash \{x\}$ such that the interior of the center curve $L$ with endpoints  $\{x,z\}$ does not intersect $\cF^{ss}(x)$.
By the transversality and completeness of the foliation $\tilde{\cF}^{cs}$, one has that for any point $w\in \cF^{cs}(x)\backslash \cF^{ss}(x)$, the strong stable leaf $\cF^{ss}(w)$ intersects $L$ in a  unique point.  Since $z\in\cF^{ss}(x)$ and  $f$ has neutral behavior along $E^c$,  we have that a subsequence of $f^n(L)$ tends  to a closed center leaf $C$ in the $C^1$-topology.  By the uniform transversality between $E^{s}$ and $E^c\oplus E^u$,  and the compactness of $M$, there exist two small positive numbers    $\epsilon$ and $\delta$  such that  for any two points $x_1,x_2$ satisfying $\ud(x_1,x_2)<\delta$,  one has that   $\cF^{ss}_{\epsilon}(x_1)$   intersects  $\cF^{cu}_{\epsilon}(x_2)$ in a  unique point which is strictly contained in
 $\cF^{ss}_{\epsilon/2}(x_1) \cap\cF^{cu}_{\epsilon/2}(x_2)$.

 We take $n$ large enough such that $f^n(L)$ is $\delta$ close to $C$, then $\cF^{ss}_{\epsilon}(f^n(L))$ intersects the annulus or M$\ddot{o}$bius band $\cF^{cu}_{\epsilon}(C)$ in a compact center curve  without boundary in the interior of $\cF^{cu}_{\epsilon}(C)$ (see Figure 1 below), which, therefore, is a compact center leaf in $f^n(\cF^{cs}(x))$, a contradiction. This ends the proof of Lemma~\ref{l.plane 3.3}.

 \begin{figure}[h]
\begin{center}
\def\svgwidth{0.3\columnwidth}
  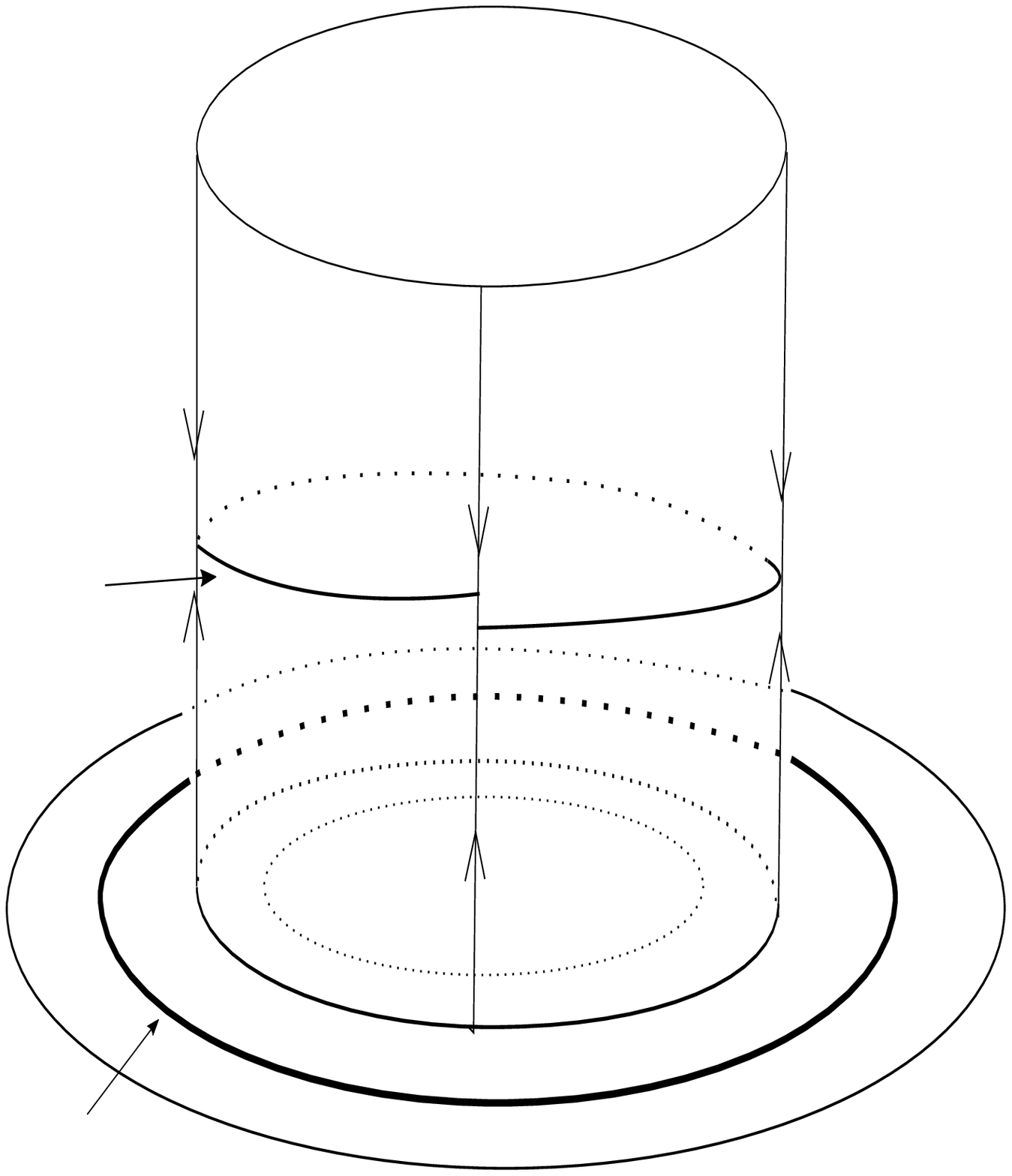
  \caption{ }
\end{center}
\end{figure}
\endproof
Now, the proof of Theorem~\ref{thm.plane and cylinder} is completed.
\endproof

At the end, we prove a property for the lifted foliations of  partially hyperbolic diffeomorphisms, which  will be used in the next section.
\begin{lemma}\label{l.unique intersection} Let $f$ be a partially hyperbolic diffeomorphism on a $3$-manifold $M$ such that $f$ has neutral behavior along the center. We denote by $\tilde{M}$ the universal cover of $M$
and by $\tilde{\cF^{i}}$ the lift of $\cF^i$ for $i=ss,cs,c,cu,uu$.

 Then for any $x,y\in \tilde{M}$, the leaf  $\tilde{\cF}^{cs}(x)$ intersects $\tilde{\cF}^{cu}(y)$ in at most one center leaf.
\end{lemma}
\proof
Assume that there exist two points $x,y\in \tilde{M}$ such that  the intersection of   $\tilde{\cF}^{cs}(x)$ and  $\tilde{\cF}^{cu}(y)$ contains two different center leaves $L_1, L_2$.

By the completeness, there exists a strong stable segment $\sigma$ whose endpoints are contained in $L_1, L_2$ respectively. Since $L_1,L_2$ are contained in the  same center unstable leaf, using a classical argument, one has that $\tilde{\cF}^{cu}$ admits a closed transversal which implies that $\cF^{cu}$ admits a null-homotopy closed  transversal.  By Theorem~\ref{thm.Novikov}, one gets that  center unstable foliation $\cF^{cu}$  has  compact leaves, contradicting to Theorem~\ref{thm.no compact leaf}.
\endproof

\section{Center flow carried by the partially hyperbolic diffeomorphisms derived from Dehn surgery: Proof of Theorem
 ~\ref{thm.conjugate to Anosov flow}}\label{construction of BPP}
In this section, we first study the properties of the partially hyperbolic diffeomorphisms in the assumption of Theorem~\ref{thm.conjugate to Anosov flow}, then we  give the proof of Theorem~\ref{thm.conjugate to Anosov flow}. At last,  we recall  the precise example in ~\cite{BPP} and give the proof of Proposition~\ref{thm.conjugation}.
\subsection{Partially hyperbolic diffeomorphisms from the Dehn surgery of Anosov flows}
Let $\phi_t$ be a non-transitive Anosov flow on an orientable  3-manifold $M$ and $\cL:M\mapsto\mathbb{R}$ be a smooth Lyapunov function of $\phi_t$.

Consider a family of transverse tori $\{T_1,\cdots,T_k\}$ which is contained in a wandering regular level $\cL^{-1}(c)$.
By the definition of a smooth Lyapunov function, one has that the regions $M^+=\cL^{-1}([c,+\infty])$ and $M^-=\cL^{-1}(-\infty,c])$ are the repelling and attracting regions for the flow $\phi_t$ respectively. We denote by $\cR$ and $\cA$ the maximal invariant sets in $M^+$ and $M^-$ respectively. By definition, one has that $\partial{M^+}=\partial{M^-}=\cL^{-1}(c)$.
We denote by $TM=E^s_X\oplus<X>\oplus E^u_X$ the hyperbolic splitting for the Anosov flow $\phi_t$.
In the following text, we denote by $W^{ss}$, $W^{cs}$, $W^c$, $W^{cu}$ and $W^{uu}$ the invariant foliations of $\phi_t$ tangent to $E_X^s$, $E^s_X\oplus<X>$, $<X>$, $<X>\oplus E^u_X$ and $E^u_X$ respectively.
\begin{proposition}\label{p.property}
Assume that  there exist $\tau>0$ and Dehn twist $\psi_i$ whose support is in $\{\phi_t(T_i)\}_{t\in(0,\tau)}$ such that the diffeomorphism $f=\psi_1\circ\cdots\psi_k\circ\phi_{\tau}$ is partially hyperbolic with one-dimensional center.  Let $TM=E^s\oplus E^c\oplus E^u$ be the partially hyperbolic splitting for $f$.

Then one has that
\begin{itemize}
\item the regions $M^+$ and $M^-$ are the repelling and attracting regions for $f$ respectively;
\item the maximal invariant sets of $f$ in $M^+$ and $M^-$ are $\cR$ and $\cA$ respectively;
\item in the region $M^+$, the bundles $<X>\oplus E^u_X$ and $E^u_X$ coincide with $E^c\oplus E^u$ and $E^u$ respectively;
\item  in the region $M^-$, the bundles $E^s_X\oplus<X>$ and $E^s_X$ coincide with $E^s\oplus E^c$ and $E^s$ respectively;
\item the splitting $E^s_X\oplus<X>\oplus E^u_X$ coincides with $E^s\oplus E^c\oplus E^u$ restricted to $\cA\cup\cR$.
\end{itemize}
\end{proposition}
\proof
Since $\psi_i$ is supported on $\{\phi_t(T_i)\}_{t\in[0,\tau]}$, one has that $f|_{M^-}=\phi_\tau$ and $f^{-1}|_{M^+}=\phi_{-\tau}$. This proves the first and the second items.

Remember that $TM=E^s_X\oplus<X>\oplus E^u_X$ is also the partially hyperbolic splitting for $\phi_{\tau}$. Since $\dim(M)=3$, by the uniqueness of partially hyperbolic splitting, one has that restricted to $\cA\cup\cR$, the splitting $E^s_X\oplus<X>\oplus E^u_X$ coincides with $E^s\oplus E^c\oplus E^u$. This proves the last item.

The following classical result ends the proof of Proposition~\ref{p.property}.
\begin{lemma} Let $U$ be an attracting region for a diffeomorphism $g$. Assume that the maximal invariant set $\La$ of $g$ in $U$ admits a dominated splitting $T_{\La}M=E\oplus F$, where $\dim(E(x))$ is a constant.  Then there exists a unique $Dg$-invariant continuous bundle $\tilde{E}$ defined on $U$ such that $\tilde{E}|_{\La}=E$.
\end{lemma}

\endproof

As a corollary, one has that
\begin{corollary} \label{c.neutral}
Let $f$ be a partially hyperbolic diffeomorphism as in the assumption of Proposition~\ref{p.property}, then $f$ has one dimensional neutral center.
\end{corollary}
The proof of Corollary~\ref{c.neutral} would be same as the proof of ~\cite[Lemma 9.1]{BPP}. As it is short,  we add the proof.
\proof  By the forth item in  Proposition~\ref{p.property}, on the set $M^-$, the bundles $E^{s}$ and $E^s\oplus E^c$ coincide with the bundles $E^s_X$ and $E^s_X\oplus <X>$. As $E^c$ is uniformly  transverse to $E^s$, in the set $M^-$,  each unit vector $v$ in $E^c$ has uniform component in the bundle $<X>$ which is bounded from above and below. Since $f$ coincides with $\phi_\tau$ on $M^-$  the diffeomorphism $f$,  the forward iterations of $v$ are uniformly  bounded from above and below.

On the set $M^+$, for the backward iterations of each unit vector, one has the analogous property.  Since the one dimensional center bundle is $f$ invariant and each orbit intersects the interior of $\cup_{i=1}^k \{\phi_t(T_i)\}_{t\in[0,\tau]}$ at most once, one gets neutral property for the center bundle.
\endproof

\begin{Remark}\label{r.center coincide}
\begin{enumerate}
\item
As the center direction  of $f$ is  the intersection of the center stable and center unstable bundles, one has that the center bundle of $f$ coincides with $<X>$ in a neighborhood of the boundary of $\{\phi_t(T_i)\}_{t\in[0,\tau]}$;
\item
Since each Dehn twist $\psi_i$ coincides with identity map in a neighborhood of the boundary of $\{\phi_t(T_i)\}_{t\in[0,\tau]}$,
by the third and fourth items of Proposition~\ref{p.property}, one has that the center stable and center unstable leaves of $f$ are invariant under $f$.
\end{enumerate}
\end{Remark}
Similar to the situation of Anosov flow $\phi_t$, one has the following corollary:
\begin{corollary}\label{co} For any $z\in M\backslash(\cA\cup\cR)$, the center leaf $\cF^c(z)$ intersects $\cL^{-1}(c)$ in a unique point.
\end{corollary}
\proof
Since the orbits of the transverse tori in $\cL^{-1}(c)$ are pairwise disjoint and the center foliation of $f$ is obtained as the intersection of the center stable and center unstable foliations, one has that each center leaf $\cF^c(z)$ intersects at most one connected component of $\cL^{-1}(c)$.

For $z\in M\backslash(\cA\cup\cR)$, one has that there exists an integer $n$ such that $f^n(z)\in\{\phi_t(\cL^{-1}(c))\}_{t\in[0,\tau]}$.
The center foliation restricted in $\{\phi_t(\cL^{-1}(c))\}_{t\in[0,\tau]}$ is obtained as the intersection of the foliations $\cF^{cs}$ and $\cF^{cu}$.
By Proposition~\ref{p.property}, one has that $$\cF^{cu}|_{\{\phi_t(\cL^{-1}(c))\}_{t\in[0,\tau]}}=f(W^{cu}|_{\{\phi_t(\cL^{-1}(c))\}_{t\in[-\tau,0]}})$$
and
$$\cF^{cs}|_{\{\phi_t(\cL^{-1}(c))\}_{t\in[0,\tau]}}=W^{cs}|_{\{\phi_t(\cL^{-1}(c))\}_{t\in[0,\tau]}}.$$
Hence, for  any $w\in\{\phi_t(\cL^{-1}(c))\}_{t\in[0,\tau]}$, the restricted center leaf $\cF^c(w)|_{\{\phi_t(\cL^{-1}(c))\}_{t\in[0,\tau]}}$ has two endpoints belonging to $T_j$ and $\phi_{\tau}(T_j)$ respectively, where $T_j$ is a connected component of $\cL^{-1}(c)$.
 Combining with the fact that $f$ equals $\phi_{\tau}$ in a neighborhood of $\{\phi_t(\cL^{-1}(c))\}_{t\in[0,\tau]}$, one has that the center leaf $\cF^c(z)$ intersects $T_i$ and $\phi_{\tau}(T_i)$ into a unique point respectively, where $T_i$ is a connected component of $\cL^{-1}(c)$; since $\phi_t(T_i)\cap T_i=\emptyset$ for $t\neq0$, one has  that the intersection of $\cL^{-1}(c)$ and $\cF^c(z)$ is a unique point.
\endproof

For the partially hyperbolic diffeomorphisms satisfying the hypothesis of Proposition~\ref{p.property}, we will show that the center foliations of such diffeomorphisms are orientable and therefore give  continuous flows (recall that the center bundle is uniquely integrable).
\begin{lemma}\label{l.center flowsss}
Let $f$ be a diffeomorphism satisfying the assumption of Proposition~\ref{p.property}.  We denote by $\cF^c$ the center foliation of $f$.
Then there exists a continuous flow $\theta_t$ on $M$ such that
\begin{itemize}
\item for any $x\in M$, one has $\Orb(x,\theta_t)=\cF^c(x)$.
\item the direction of flow $\theta_t$ gives the same transverse orientation to $\cL^{-1}(c)$ as the direction of flow $\phi_t$ does.
\item restricted to $\cA\cup\cR$, the orientation of $E^c$ given by the direction of flow $\theta_t$ coincides with the orientation given by the direction of the flow $\phi_t$.
\end{itemize}
\end{lemma}
\proof By Remark ~\ref{r.center coincide}, one has that the center foliation $\cF^c$ coincides with the center foliation $W^c$ in a neighborhood of $\cL^{-1}(c)$.

For each center leaf intersecting $\cL^{-1}(c)$, we give it the same orientation as the one of $W^c$ given by the flow direction. Since in the region $M^-$, the center stable  and strong stable foliations of $f$ coincide with the ones of the Anosov flow $\phi_t$, by the fact that $\cF^c$ is  transverse to the strong stable direction in each center stable leaf, one has that in the region $M^-\backslash \cA$, restricted to each center stable leaf, each  leaf of $\cF^c$  cuts the strong stable leaf with the same orientation as $W^c$ does; another way to observe this is to lift them to the universal cover, and one still has that  the lift of $\cF^c$ coincides with the lift of $W^c$ in a neighborhood of the lift of $\cL^{-1}(c)$ and the lift of strong stable foliations $\cF^{ss}$, $W^{ss}$ restricted to the lift of  $M^-\backslash \cA$ also coincide.  Combining with Corollary~\ref{co}, one has that the orientation of  the $\cF^c$-leaves intersecting $\cL^{-1}(c)$ induces the same orientation of the center leaves in $\cA$  as the the one given by the flow. One can apply the same argument to the region $M^+$. Finally, one gets that the center foliation $\cF^c$ is orientable,  and one can give it an orientation such that
\begin{itemize}
\item it gives the same orientation on the set $\cA\cup\cR$ as the one given by the flow $\phi_t$;
\item it gives the same transverse  orientation to $\cL^{-1}(c)$  as the one given by the flow $\phi_t$.
\end{itemize}
Hence, the center foliation $\cF^c$ can give a continuous flow, denoted as $\theta_t$ and satisfying the posited properties. This ends the proof of Lemma~\ref{l.center flowsss}.
\endproof

Now, we are ready to give the proof of Theorem~\ref{thm.conjugate to Anosov flow}.

\proof[Proof of Theorem~\ref{thm.conjugate to Anosov flow}]
By Corollary ~\ref{c.neutral}, the partially hyperbolic diffeomorphism $f$ has neutral center.

Let $\theta_t$ be the continuous flow given by
 Lemma~\ref{l.center flowsss} with respect to the center foliation $\cF^c$ of $f$.

We shall build the conjugation between the center flow $\theta_t$ and the original Anosov flow $\phi_t$.
 We shall first define the conjugation on the set $M^-\backslash \cA$ which is Id restricted to the boundary of $M^-$, then we extend
  it to be a homeomorphism of $M^-$ which equals Id on $\cA$. Similarly, we define the conjugation on the set $M^+$. In the end, we get the conjugation between $\theta_t$ and $\phi_t$.

 We denote  by $\cF^{ss}$, $\cF^{cs}$, $\cF^{cu}$ and $\cF^{uu}$ the strong stable, center stable, center unstable and strong unstable foliations of $f$ respectively.

Let  $\pi:\tilde{M}\mapsto M$ be  the universal cover of $M$. We denote  by $\tilde{\cF}^l$ and $\tilde{W}^l$ the
lifts   of the foliations $\cF^l$ and  $W^l$ respectively, for any $l=ss,cs,c,cu,uu$. Given a foliation $\cF$ on $M$ and a submanifold $M^\prime\subset M$,
the leaf of $\cF|_{M^\prime}$ through a point $x\in M^\prime $ is the connected component of $\cF(x)\cap M^{\prime}$ containing $x$.

 \begin{proposition}\label{p.conjugation} There exists a homeomorphism $h^s:M^-\mapsto M^-$ such that
 \begin{itemize}
  \item the map $h^s$ preserves every leaf of the   foliation $W^{cs}|_{M^-}=\cF^{cs}|_{M^-}$;
 \item the map $h^s$ takes the orbits of the flow $\phi_t|_{M^-}$ to the ones of $\theta_t|_{M^-}$;
 \item the map $h^s$ coincides with $\Id$ on the set $\cA\cup \cL^{-1}(c)$.
 \end{itemize}
 \end{proposition}

 From now on, for each set $A\subset M$, we denote by $\tilde{A}=\pi^{-1}(A)$.
To prove Proposition~\ref{p.conjugation}, we need the following lemma:
\begin{lemma}\label{l.equivalent intersection}
 For any points $x\in\tilde{M}^{+}$ and $y\in\tilde{M}^{-}$, the center unstable leaf $\tilde{W}^{cu}(x)$ intersects
the center stable $\tilde{W}^{cs}(y)$ non-empty if and only if the center unstable leaf $\tilde{\cF}^{cu}(x)$ intersects
the center stable $\tilde{\cF}^{cs}(y)$ non-empty. More precisely, we have the following:
$$z\in\tilde{W}^{cu}(x)\cap \tilde{W}^{cs}(y)\cap\widetilde{\cL^{-1}(c)}\Longleftrightarrow z\in\tilde{\cF}^{cu}(x)\cap \tilde{\cF}^{cs}(y)\cap\widetilde{\cL^{-1}(c)}.$$
\end{lemma}

\proof[Proof of Lemma~\ref{l.equivalent intersection}]
First, we need the following result:

\begin{claim}\label{c.equivalent intersection} Given two points  $p,q\in M^-$, if $p\in\cF^{cs}(q)$, then $p,q$ are  in the same connected component of $\cF^{cs}(q)\cap M^-$.
\end{claim}
\proof By the fourth item of Proposition~\ref{p.property}, one has that $\Orb^+(p,\phi_t)\cup\Orb^+(q,\phi_t)$ is contained in $\cF^{cs}(q)$.
Since the orbits of $p,q$ under the flow $\phi_t$ converge to an orbit in $\cA$, once again by the fourth item of Proposition ~\ref{p.property}, one has that $p,q$ are  in the same connected component of $\cF^{cs}(q)\cap M^-$.
\endproof

Given $x\in\tilde{M}^{-}$ and $y\in\tilde{M}^{+}$, by Lemma~\ref{l.unique intersection},
 the intersection of $\tilde{W}^{cs}(x)$ and $\tilde{W}^{cu}(y)$ consists of  at most one leaf  of $\tilde{W}^c$, and the same property holds for the foliations $\tilde{\cF}^{cs}$ and $\tilde{\cF}^{cu}$.
If the intersection of  $\tilde{W}^{cs}(x)$ and $\tilde{W}^{cu}(y)$ is not empty, by the choices of $x,y$,  there exists a  unique  point $z\in\widetilde{\cL^{-1}(c)}$  contained in $\tilde{W}^{cu}(x)$ and $\tilde{W}^{cs}(y)$,
which is equivalent to that $z\in\widetilde{\cL^{-1}(c)}$  is contained in $\tilde{\cF}^{cu}(x)$ and $\tilde{\cF}^{cs}(y)$; this is due to Claim~\ref{c.equivalent intersection} and the fact that   these foliations coincide on the corresponding region. This ends the proof of Lemma~\ref{l.equivalent intersection}.
\endproof

Now, we are ready to  give the proof of Proposition~\ref{p.conjugation}.
\proof[Proof of Proposition~\ref{p.conjugation}]
Given  $x\in\tilde{M}^{-}\backslash\tilde{\cA}$, there exists a point $y\in\tilde{M}^{+}\backslash\tilde{\cR}$ such that $x$ is on the
 positive orbit of $y$ for the lifted flow  $\tilde{\phi}_t$; by Corollary~\ref{co}, one has the analogous property for the flow $\tilde{\theta}_t$. Let $z$ be the unique intersection between $\Orb(x,\tilde{\phi}_t)=\tilde{W}^c(x)$ and the transverse section $\pi^{-1}(\cL^{-1}(c))$, then  one has that
  $$z\in\tilde{W}^{cu}(y)\cap \tilde{W}^{cs}(x)\cap\widetilde{\cL^{-1}(c)}.$$
  By Lemma~\ref{l.equivalent intersection}, one has that
$$z\in\tilde{\cF}^{cu}(y)\cap \tilde{\cF}^{cs}(x)\cap\widetilde{\cL^{-1}(c)}.$$
Since each leaf of the foliation $\tilde{\cF}^{cs}$ is a plane, by  completeness and  transversality,  the center   leaf $\tilde{\cF}^c(z)$  intersects the strong stable leaf $\tilde{\cF}^{ss}(x)$ in a  unique point and we denote it as $h^s(x)$ (as it is shown in Figure 2).
Similarly, we can define a map $\tau^s$ by exchanging the roles  of $(\tilde{\phi_t}, \tilde{W}^{ss})$ and $(\tilde{\theta}_t,\tilde{\cF}^{ss})$ in the definition of $h^s(x)$.
\begin{figure}[h]
\begin{center}
\def\svgwidth{0.5\columnwidth}
  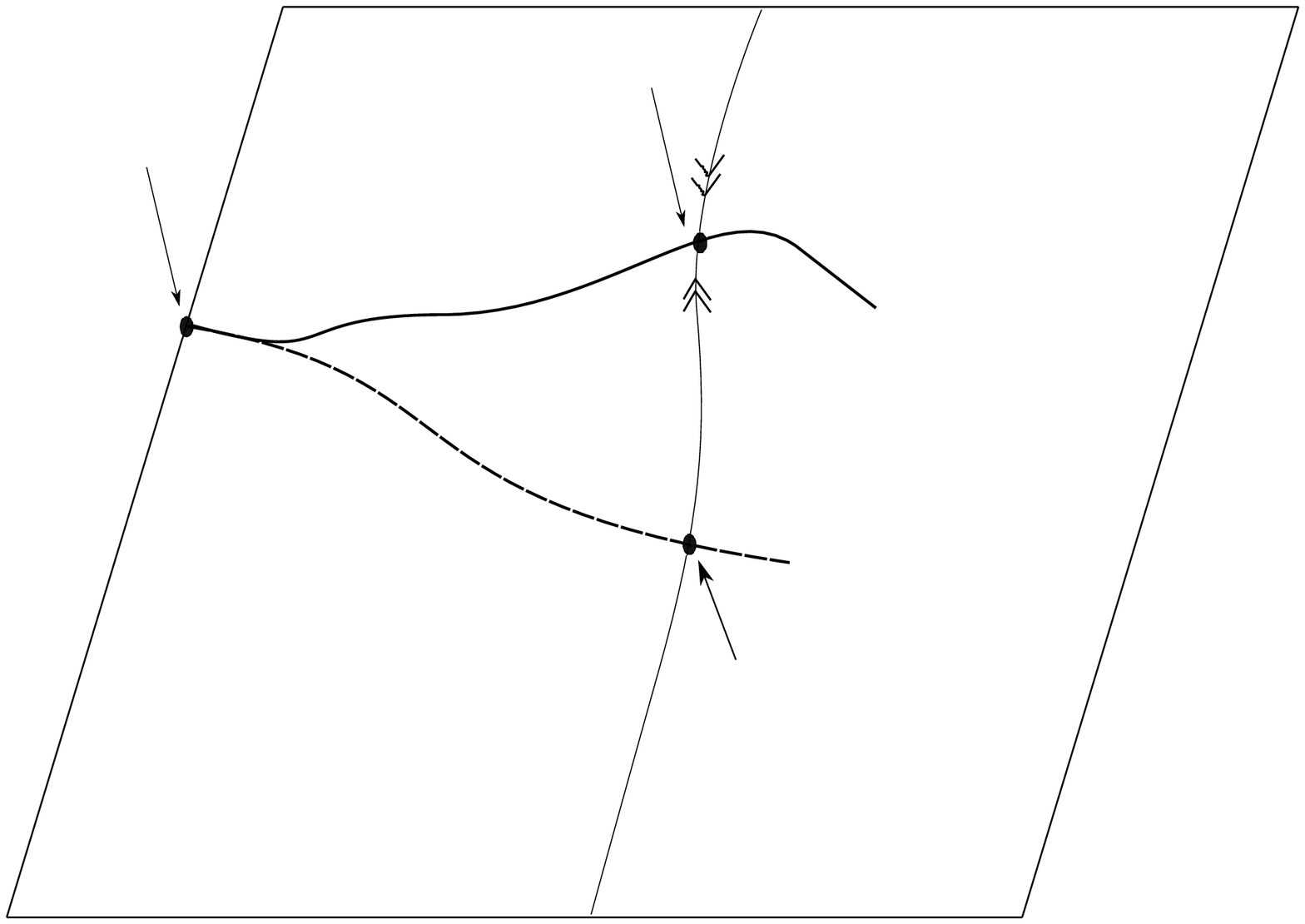
  \caption{ }
\end{center}
\end{figure}

By   definition, one can check that the maps $h^s,\tau^s : \tilde{M}^-\backslash \tilde{\cA}\mapsto \tilde{M}\backslash(\tilde{\cA}\cup\tilde{\cR})$ are continuous  and commutative with the automorphisms of $\tilde{M}$ induced by $\pi_1(M)$. Moreover,  one has that
\begin{itemize}
\item the maps $h^s, \tau^s$ coincide with $\Id$ in a neighborhood of  $\widetilde{\cL^{-1}(c)}$ restricted in $\tilde{M}^-$;
\item the maps $h^s, \tau^s$ preserve the orientation of the foliations $\tilde{\cF}^{cs}|_{\tilde{M}^-}, \tilde{\cF}^{ss}|_{\tilde{M}^-}$ and send the positive orbits of one  center flow to the positive orbits of  the other.
\end{itemize}
\begin{claim}\label{c.injective} The maps $h^s$ and $\tau^s$ are injective.
\end{claim}
\proof Assume that $h^s$ is not injective, then there exist two different points $x,y\in\tilde{M}^-\backslash\tilde\cA$ such that $h^s(x)=h^s(y)$.
By definition, one has that $x\in\tilde{\cF}^{ss}(y)$. Moreover, $x$ and $y$ are on the same orbit of the Anosov flow $\tilde{\phi}_t$. By the forth item in Proposition~\ref{p.property},  the orbit segment from $x$ to $y$ for the flow $\tilde{\phi}_t$ is contained in the leaf $\tilde{\cF}^{cs}(x)$ and is transverse to $\tilde{\cF}^{ss}$. Since the leaf $\tilde{\cF}^{cs}(x)$ is a plane and the leaf $\tilde{\cF}^{ss}(x)$ intersects $L(x,y)$ twice, we get the contradiction.

The same argument applies for $\tau^s$.
\endproof

\begin{claim}\label{c.conjugation} The images of $h^s$ and $\tau^s$ are contained in  $\tilde{M}^-\backslash\tilde{\cA}$.
\end{claim}
\proof  Let us first recall that
$$\tilde{\cF}^{cs}|_{\tilde{M}^{-}}=\tilde{W}^{cs}|_{\tilde{M}^-} \textrm{ and }
\tilde{\cF}^{cu}|_{\tilde{M}^{+}}=\tilde{W}^{cu}|_{\tilde{M}^+}.$$
By definition, for any point $x\in\tilde{M}^-\backslash\tilde{\cA}$, $h^s$ maps the connected component of $\tilde{\cF}^{cs}(x)\cap (\tilde{M}^-\backslash\tilde{\cA})$ containing $x$ to a connected component of $\tilde{\cF^{cs}}(x)\cap (\tilde{M}\backslash\tilde{\cA})$. Notice that each connected  component of $\tilde{\cF^{cs}}(x)\cap (\tilde{M}\backslash\tilde{\cA})$ has non-empty interior  in $\tilde{\cF}^{cs}(x)$ and its boundary consists of  center leaves  in $\tilde{\cA}$. Since $h^s$ coincides with identity in a neighborhood of $\widetilde{\cL^{-1}(c)}$ restricted in $\tilde{M}^-$, by the continuous and injective property of $h^s$,
one has that $h^s$ maps the connected component of $\tilde{\cF}^{cs}(x)\cap (\tilde{M}^-\backslash\tilde{\cA})$ containing $x$ into itself,  for any $x\in\tilde{M}^-\backslash\tilde{\cA}$.

One can prove the claim for $\tau^s$ analogously, ending the proof of Claim~\ref{c.conjugation}.
\endproof

\begin{claim}~\label{c.same strong stable}
For any points $x, y\in\tilde{M}^-$, one has that
$$y\in\tilde{\cF}^{ss}(x)\Longleftrightarrow y\in\tilde{W}^{ss}(x).$$
\end{claim}
\proof Since  $f|_{M^-}=\phi_{\tau}|_{M^-}$, one can take a lift $\tilde{f}$ of $f$ such that the forward orbit of $x$
under $\tilde{f}$ coincides with the one under $\tilde{\phi}_{\tau}$. If $y\in\tilde{\cF}^{ss}(x)\cap \tilde{M}^-$,
then the forward orbit of $y$ under $\tilde{f}$ coincides with the one under $\tilde{\phi}_{\tau}$ since
 $f|_{M^-}=\phi_{\tau}|_{M^-}$. Hence,
the distance $\ud(\tilde{f}^n(x),\tilde{f}^n(y))=\ud(\tilde{\phi}_{n\tau}(x),\tilde{\phi}_{n\tau}(y))$ tends to zero exponentially. Since $\tilde{\phi}_\tau$ has neutral behavior along the center, one has that  $y\in\tilde{W}^{ss}(x)$.  One can argue for the other side  analogously, concluding.
\endproof
 By the definitions of $h^s$ and $\tau^s$, one has that $h^s(x)\in \tilde{W}^{ss}(x)$ and $\tau^s(x)\in\tilde{\cF}^{ss}(x)$.  By Claims~\ref{c.conjugation},~\ref{c.same strong stable} and the definitions of $h^s,\tau^s$, one has that $\tau^s\circ h^s=h^s\circ\tau^s=\Id$. Hence, $h^s$ is an orientation preserving homeomorphism on the set $\tilde{M}^-\backslash\tilde{A}$.  Moreover, $h^s$ maps each leaf of $\tilde{\cF}^{ss}|_{\tilde{M}^-\backslash\tilde{\cA}}$ into itself surjectively as a homeomorphism, hence one can extend  $h^s$  to $\tilde{\cA}$ as $\Id$ such that  $h^s$ is continuous along each leaf of $\tilde{\cF}^{ss}|_{\tilde{M}^-}$.
  \begin{lemma} The maps $h^s,\tau^s: \tilde{M}^-\mapsto \tilde{M}^-$ are continuous.
  \end{lemma}
  \proof One only needs to prove that the map $h^s$ is continuous at $\tilde{\cA}$.  The case for $\tau^s$ follows analogously.

  Assume, on the contrary, there exists a point $x_0\in\tilde{\cA}$
  where $h^s$ is not continuous.  Then there exist $\epsilon_0>0$
   and a sequence of points $\{x_n\}_{n>0}\subset\tilde{M}^-\backslash\tilde{\cA}$ such that
 $$\lim_{n\rightarrow\infty} x_n=x_0 \textrm{ and } \ud(h^s(x_n),h^s(x_0))>\epsilon_0.$$

 To continue the proof, we need the following result:
 \begin{claim}~\label{c.123}
 On the universal cover $\tilde{M}$, for the lifts of $f$,  each center unstable leaf intersects a strong stable leaf in at most one point. The analogous property also holds for $\phi_\tau$.
 \end{claim}
 \proof If there exists a strong stable leaf intersecting a center unstable leaf in two points,
  by a classical argument, one gets a closed transversal for the center unstable foliation, which implies that  the center unstable foliation for the diffeomorphism on $M$ admits a null-homotopy closed transversal. By Novikov's theorem, the center unstable foliation for the diffeomorphism on $M$ has compact  leaves, contradicting to Theorem~\ref{thm.no compact leaf}.
 \endproof
 Since $x_n$ tends to $x_0$, by Claim~\ref{c.123}, the center unstable leaf $\tilde{W}^{cu}(x_n)$ intersects the strong
 stable leaf $\tilde{W}^{ss}_{\epsilon_0}(x_0)=\tilde{\cF}_{\epsilon_0}^{ss}(x_0)$ in a unique point $y_n$, for $n$ large.
  Then $y_n$ tends to $x_0$. Moreover, since each connected component of $\tilde{W}^{cu}(x_n)\cap\tilde{M}^-$ is the forward $\tilde{\phi}_t$-orbit of a
 connected component of $\tilde{W}^{cu}(x_n)\cap \widetilde{\partial{M}^-}$,
 the backward orbits of $x_n$ and $y_n$ under the flow $\tilde{\phi}_t$ intersect the same connected component $P_n$ of $\widetilde{\partial{M}^-}\cap \tilde{W}^{cu}(x_n)$ into unique points, and we denote them by  $p_n$ and $q_n$ respectively, where $P_n$ is diffeomorphic to $\mathbb{R}$.
 Since $\tilde{W}^{cu}$ coincides with $\tilde{\cF}^{cu}$ in a neighborhood of $\widetilde{\partial{M}^-}$, one has that the forward orbits of $p_n$ and $q_n$ under the flow $\tilde{\theta}_t$ are on the same center unstable leaf $\tilde{\cF}^{cu}(q_n)$. By the choices of $p_n, q_n$ and Claim~\ref{c.123}, one has that the center unstable leaf $\tilde{\cF}^{cu}(q_n)$ intersects $\tilde{\cF}^{ss}(x_0)$ and $\tilde{\cF}^{ss}(x_n)$ into  unique points; remember that the orbits of the flow $\tilde{\theta}_t$ also lie in $\tilde{\cF}^{cu}$, by the definition of $h^s$, one has that their intersections must be  $h^s(y_n)$ and $h^s(x_n)$ respectively.
By definition, $h^s$ is continuous restricted to the strong stable leaf $\tilde{\cF}^{ss}(x_0)$, hence $\ud(h^s(y_n),x_0)$ tends to zero. By the continuity of center unstable foliation $\tilde{\cF}^{cu}$, the center unstable plaque $\tilde{\cF}_{\epsilon_0/3}^{cu}(h^s(y_n))\subset\tilde{\cF}^{cu}(q_n)$ intersects the strong stable leaf $\tilde{\cF}^{ss}(x_n)$ into a point $z_n$ for $n$ large, hence
 $$\ud(z_n, x_0)<\ud(z_n, h^s(y_n))+\ud(h^s(y_n),x_0)<\epsilon_0, \textrm{  for $n$ large}.$$
  Once again,  by Claim ~\ref{c.123}, one has that $z_n=h^s(x_n)$,  contradicting to $\ud(h^s(x_n), x_0)<\epsilon_0$.
  \endproof

  Now, the map $h^s:\tilde{M}^-\mapsto\tilde{M}^-$ is a homeomorphism.
  By the  definition of $h^s$,   one
   \vspace{1.5mm}
  \\ has that $h^s$ maps the orbits of $\tilde{\phi}_t|_{\tilde{M}^-}$ to the orbits of $\tilde{\theta}_t|_{\tilde{M}^-}$, and preserves the orientation of the orbits. Since $h^s$ commutes with the automorphisms  on $\tilde{M}$ induced by $\pi_1(M)$, the projection of $h^s$ on the base manifold defines a homeomorphism of $M^-$ satisfying the announced properties, ending the proof Proposition~\ref{p.conjugation}.
\endproof

\proof[Ending the proof of Theorem~\ref{thm.conjugate to Anosov flow}]
By applying Proposition~\ref{p.conjugation} to the reversed dynamics on the set $M^+$, one gets a homeomorphism $h^u:M^+\mapsto M^+$ satisfying the analogous properties.
We define a homeomorphism $h: M\mapsto M$ in the following way:
 \begin{displaymath}
 h(x)=\left\{\begin{array}{ll}
 h^s(x) &\textrm{ if $x\in M^-$}
 \\
 h^u(x) &\textrm{ if $x\in M^+$}
 \end{array}\right.
\end{displaymath}
The homeomorphism $h$ coincides with $\Id$ on the set $\cA\cup \cR\cup \cL^{-1}(c)$.  One can check that $h$ sends the orbits of $\phi_t$
to the orbits of $\theta_t$ and preserves the orientation of the flows.
This  proves that $\theta_t$ is topologically equivalent to Anosov flow $\phi_t$.
\endproof

\section{The anomalous example in ~\cite{BPP}:  proof of Proposition~\ref{thm.conjugation}}
\subsection{Construction of  the  example in ~\cite{BPP}}

  In ~\cite[Section 4]{BPP}, the authors  built a 3-manifold $N$ supporting  a smooth non-transitive Anosov flow $\psi_t$ having two transverse tori $T_1$ and $T_2$ with the properties:
  \begin{enumerate}
  \item[(P1).] $T_1\cup T_2$ is far away from the non-wandering set of $\psi_t$;
  \item[(P2).] $T_1\cup T_2$ separates $N$ into two connected components $N^+$ and $N^-$;
  \item[(P3).] $N^+$ is a repelling region of $\psi_t$ and the maximal invariant set of $\psi_t$ in $N^+$ is a repeller $\cR$;
  \item[(P4).]    $N^-$ is an attracting region of $\psi_t$ and the maximal invariant set of $\psi_t$ in $N^-$ is an attractor $\cA$.
  \end{enumerate}
  By ~\cite{HP}, one has the the stable and unstable foliations of $\psi_t$ are $C^1$, hence on the transverse torus $T_i$, the stable manifold of $\psi_t$ and the unstable manifold of $\psi_t$ induce two  $C^1$ foliations $\cF^s_i$ and $\cF^u_i$ respectively. For these foliations, one has the following property:
  \begin{enumerate}
  \item[(P5).]the induced foliations $\cF^s_i$ and $\cF^u_i$ consist of two Reeb components.
  \end{enumerate}
   The compact leaves of these Reeb components belong to the stable or the unstable manifolds of the periodic orbits of $\psi_t$, since these compact leaves are transverse to the flow in the stable or the unstable manifolds of the flow $\psi_t$. One can take a $C^1$ coordinate $\theta_i$ for each transverse torus  so that under such coordinates,  the induced foliations $\cF^s_i$ and $\cF^u_i$ on transverse torus $T_i$ are exactly as shown in  Figure 3 (for details see Lemma 4.1 in ~\cite{BPP}).
  \begin{figure}[h]\label{f.two reeb components}
\begin{center}
\def\svgwidth{0.25\columnwidth}
  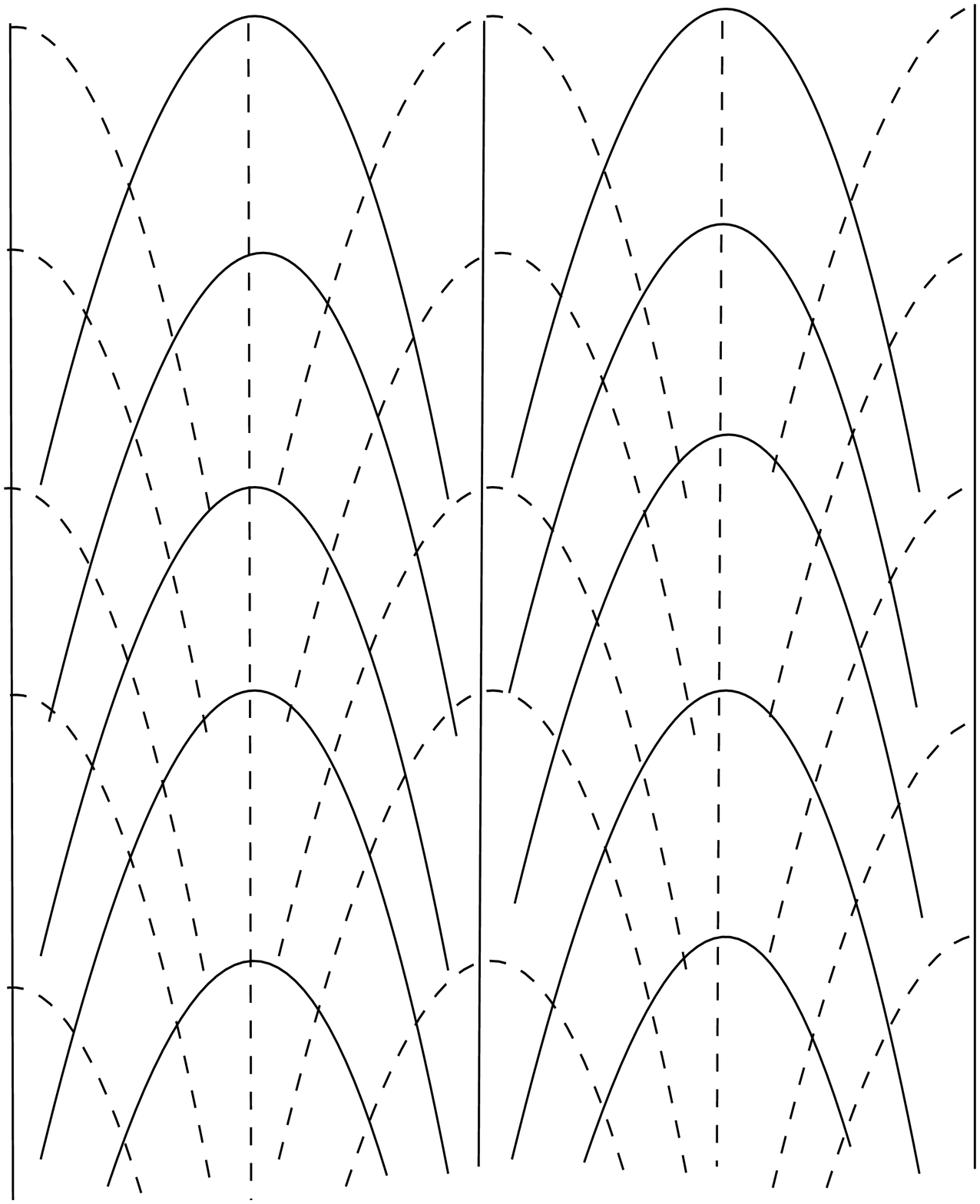
  \caption{The real lines and the dash lines  denote the leaves of the lifts of foliations $\cF^s_i$ and $\cF^u_i$  on the universal cover respectively.}
\end{center}
\end{figure}

Under this coordinate, one considers the Dehn twist on $T_1\times[0,1]$ with the form $\tilde{\Psi}:(x,t)\mapsto(x+(0,\alpha(t)),t)$, where $\alpha(t)$ is a smooth non-decreasing bump function supported on $[0,1]$  such that $\alpha(0)=0$ and $\alpha(1)=1$, and the circle $\{(0,t)\}_{t\in[0,1]}$ on $T_1$ is in the homotopy class of the compact leaves of $\cF^s_i$.  We denote by $\Phi_t(x)=x+(0,\alpha(t))$ which is a smooth diffeomorphism on $T_1$. Under this coordinate, one has that $\Phi_t(\cF^u_i)\pitchfork\cF^s_i$ for any $t\in[0,1]$. Consider the Dehn twist $\Theta_n=\Gamma_n^{-1}\circ\tilde{\Psi}\circ\Gamma_n$ defined on $\{\psi_t(T_1)\}_{t\in[0,n]}$, where $\Gamma_n:\{\psi_t(T_1)\}_{t\in[0,n]}\mapsto T_1\times[0,1]$ is of the form $(\psi_t(x))\mapsto(x,\frac{t}{n})$.
In ~\cite{BPP}, the authors prove that
\vspace{1mm}
\begin{theorem}~\cite[Theorem 8.1 and Lemma 9.1]{BPP} For $n>0$ large, the diffeomorphism $f_b=\Theta_n\circ\psi_n$ is partially hyperbolic with one dimensional neutral center.
\end{theorem}

   By ~\cite[Proposition 1.9]{BZ}, there exists a smooth Lyapunov function such that $\{T_1,T_2\}$ is a wandering regular level of this Lyapunov function.
 As a consequence, one can apply Proposition ~\ref{p.property} to the diffeomorphism $f_b$.

\subsubsection{Action of $f_b$ on the space of center leaves intersecting $T_1$}

We denote by $\cF^l$  the $f_{b}$-invariant foliation tangent to $E^l$, for $l=ss,cs,c,cu,uu$.
For each transverse torus $T_i$, we lift the foliations $\cF^s_i,\cF^u_i$ to the universal cover $\R^2$, and we denote them as $\tilde{\cF}^s_i$ and $\tilde{\cF}^{u}_i$ respectively.  By transversality, a $\tilde{\cF}^s_i$-leaf intersects a $\tilde{\cF}^u_i$-leaf  in at most one point.
Since under the coordinate of $\theta_1$, the induced foliations on $T_1$ are shown as in Figure 3, one has the following result (as it is shown in Figure 4):
    \begin{lemma}\label{l.lift intersection} Under the coordinate $\theta_1$ of $T_1$, we lift the induced foliations to the universal cover $\R^2$. Let $\mathcal{T}\in\diff^1(\R^2)$ be the translation of the form $(t,s)\mapsto (t,s+1)$. Then for any $x\in\R^2$, one has that  the leaf $\mathcal{T}( \tilde{\cF}^u_i(x))$ intersects $\tilde{\cF}^s_i(x)$  in a unique point.
    \end{lemma}

   \begin{figure}
\begin{center}
\def\svgwidth{0.4\columnwidth}
  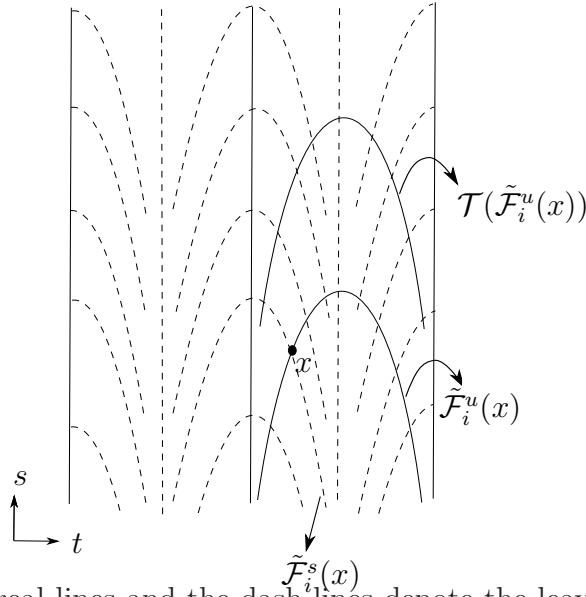
  \caption{The real lines and the dash lines  denote the leaves of the lifts of foliations $\cF^u_i$ and $\cF^s_i$  to the universal cover respectively.}
\end{center}
\end{figure}

Now, one can define a homeomorphism $\eta\in\homeo(\R^2)$ which maps the point $x$ to the point $\cT(\tilde{\cF}^u_1(x))\cap\tilde{\cF}^s_1(x)$. One can check that $\eta$ induces a homeomorphism on $T_1$ and for notational convenience, we still denote it as $\eta$. By definition, one has that $\eta$ keeps every leaf of $\cF^s_1$ and every leaf of $\cF^u_1$ invariant.
 Moreover,  the homeomorphism $\eta$ coincides with identity on union of the compact leaves of $\cF^s_1$ and the union of the compact leaves of $\cF^u_1$.

\begin{lemma}\label{l.existence of fixed center leaf}
 The action of the  diffeomorphism $f_b$  on the space of center leaves intersecting $T_1$ is equivalent to the homeomorphism $\eta^{-1}$.
\end{lemma}
\proof
By Proposition~\ref{p.property},
   in $W_1=\{\psi_t(T_1)\}_{t\in[0,n]}$, the center unstable foliation $\cF^{cu}$ is given by $\Theta_n(W^{cu})$ where $W^{cu}$ is the center unstable foliation of $\psi_n$. Consider the $C^1$ coordinate $\vartheta=(\theta_1, \Id)\circ \Gamma_n$. Under this coordinate and restricted to  $W_1$, the center unstable foliation $W^{cu}$ coincides with the product foliation $\cF^u_1\times [0,1]$.  Then,
    by Lemma~\ref{l.lift intersection}, for each point $x\in T_1$, the center leaf through $x$ intersects $\psi_n(T_1)$ into $\psi_n(\eta(x))$; one can observe this by lifting $W_1$ to the universal cover $\tilde{W}_1$; and the lifts of the non-compact leaves of $\cF^{cu}|_{W_1}$
intersect the lifts of the non-compact leaves of $\cF^{cs}|_{W_1}$ in the way  shown in Figure 5.

\begin{figure}[h]\label{f.intersection of center stable and center unstable}
\begin{center}
\def\svgwidth{0.3\columnwidth}
  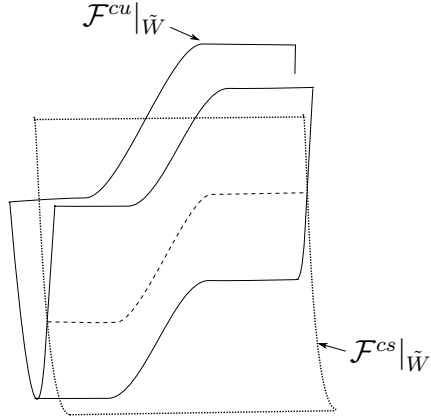
  \caption{The dash line denotes the center leaf obtained by the intersection of center stable and center unstable leaves.}
\end{center}
\end{figure}

Since  $\Theta_n$ coincides with $\Id$ in a neighborhood of the boundary of $\{\psi_{t}(T_i)\}_{t\in[0,n]}$, the diffeomorphism $f_{b}$ sends the center leaf through $x\in T_1$ to the center leaf through $\psi_{n}(x)$, and  the center leaf through the point $\psi_{n}(x)$ is the one through the point $\eta^{-1}(x)\in T_1$.
\endproof

For the transverse torus $T_2$, the diffeomorphism $f_b$ coincides with $\psi_n$ in $\{\psi_t(T_2)\}_{t\in \R}$, therefore, in $\{\psi_t(T_2)\}_{t\in \R}$, the center stable and center unstable foliations of $f_b$ coincide with the ones of $\psi_n$. Hence, the center leaves intersecting $T_2$ are  invariant under $f_b$.

As a Corollary of Lemma~\ref{l.existence of fixed center leaf}, one has the following:
\begin{corollary}   For the partially hyperbolic diffeomorphism $f_b$,  for   any $j=1,2$ and any leaf  $\cF^s_j(y)$ (resp. $\cF^u_j(y)$) of  the foliation $\cF^s_j$ (resp. $\cF^u_j$),  there exists a point $x\in\cF^s_j(y)$ (resp. $\cF^u_j(y)$)  such that $f_b$ preserves the center leaf through $x$, that is,  $f_b(\cF^c(x))=\cF^c(x)$.
\end{corollary}

\subsection{Robust completeness of the center stable foliation of ~\cite{BPP} example}
By Theorem~\ref{thm.plane and cylinder}, the center stable and center unstable foliations of $f_b$ are complete. Indeed, for this particular example, one can get the robust completeness of the center stable and center unstable foliations.

Now, we will use Theorem~\ref{thm.conjugate to Anosov flow} to give the proof of Proposition~\ref{thm.conjugation}.
\proof[Proof of Proposition~\ref{thm.conjugation}]
Since $f_b$ has neutral behavior along the center, by Theorem~\ref{r.preserves all the leaves},  there exists a $C^1$ small neighborhood $\cV$ of $f_b$ such that for any $g\in\cV$, one has that
\begin{itemize}
\item $g$ is dynamically coherent;
 \item there exists a homeomorphism $h_g:N\mapsto N$ such that for any $x\in N$ and  $i=c,cs,cu$, one has
$$h_g(\cF^i(x))=\cF^i_g(h_g(x))\textrm{ and }
           h_g(f_b(\cF^i(x)))=g(\cF_g^i(h_g(x))) ;$$
\item the homeomorphism $h_g$ tends to identity in the $C^0$-topology when $g$ tends to $f$.
     \end{itemize}

Recall that  the maximal invariant sets  of $f_b$ in $N^+$ and $N^-$ are
$\cR$ and $\cA$ respectively. By Theorem 7A.1 in ~\cite{HPS},  one can choose  a small enough neighborhood $\cU\subset\cV$ of $f_b$ such that for any $g\in\cU$,  the maximal invariant set of $g$ in the region  $N^{-}$ (resp. $N^+$) is  $h_g(\cA)$ (resp. $h_g(\cR)$).
Hence, the chain recurrent set of $g$ is contained in $h_g(\cA\cup\cR)$.

    We will first show  the following:
    \begin{lemma}\label{l.completeness in repeller} For any point $x\in h_g(\cR)$, one has that
    $$\cF^{ss}_g(\cF^c_g(x))=\cF^{cs}_g(x).$$
    \end{lemma}
    \proof
    Notice that $\cF^{cs}_g(x)=h_g(\cF^{cs}(h^{-1}_g(x)))$. We lift  the foliations to the universal cover $\tilde{N}$. Let $\tilde{\cF}_g^c(y)$
    be a lift of $\cF^c_g(x)$. Notice that  every center leaf in $h_g(\cR)$ is $g$-invariant.
     \begin{claim}\label{c.keep center leaf invariant} There exists a lift $\tilde{g}$ of $g$ such that  every center leaf in $\tilde{\cF}^{cs}_g(y)$ is $\tilde{g}$-invariant.
     \end{claim}
     \proof By leaf conjugacy, there is at most one compact center leaf  contained in $\cF_g^{cs}(x)$. We take a non-compact center leaf
      $L\subset  \cF_g^{cs}(x)$,  and we denote by $\tilde{L}\subset \tilde{\cF}^{cs}_g(y)$ a lift of $L$.  Then there exists a unique lift
     $\tilde{g}$ of $g$ such that $\tilde{L}$ is $\tilde{g}$-invariant.

     For any center leaf $\tilde{\cF}_g^c(z)\subset \tilde{\cF}^{cs}_g(y)$, by the fact that $\cF^c$ is topologically Anosov,
     there exists a strong stable curve $\tilde{\sigma}(t)_{t\in[0,1]}\subset \tilde{\cF}^{cs}_g(y)$  such that $\tilde{\sigma}(0)\in\tilde{L}$ and
     $\tilde{\sigma}(1)\in \tilde{\cF}^c_g(z)$.
     Let $\sigma(t)$ be the projection of  $\tilde{\sigma}(t)$ on $\cF^{cs}_g(x)$, by the invariance of the center leaves, one has that there exists a continuous  family of center curves $\gamma_t(s)$ joining the curve $\sigma(t)$ with $g(\sigma(t))$. Now, we lift the family of center  curves to the leaf $\tilde{\cF}^{cs}_g(y)$, then one gets  a continuous  family of center curves $\tilde{\gamma}_t(s)$ joining $\tilde{\sigma}(t)$ to a lift $\alpha(t)$ of $g(\sigma(t))$.
     By the uniqueness of $\tilde{g}$ and the non-compactness of $L$,  the lift of $g(\sigma(0))$ can only be $\tilde{g}(\tilde{\sigma}(0))$ which implies that $\alpha(t)$
     can only be  $\tilde{g}(\tilde{\sigma}(t))$. Hence, one has that $\tilde{\cF}_g^c(z)$ is invariant under $\tilde{g}$.
     \endproof

We only need to prove that,  restricted to the leaf  $\tilde{\cF}^{cs}(y)$,  every strong stable leaf  intersects every center leaf.
Assume, on the contrary, that $\tilde{\cF}^{ss}_g(\tilde{\cF}^c_g(q))$ has boundary leaves for some $q\in\tilde{\cF}^{cs}(y)$. Let $\tilde{\cF}^{ss}_g(p)$ be one of the boundary leaf.
\begin{claim} Let $\tilde{g}$ be the lift of $g$ given by Claim~\ref{c.keep center leaf invariant}, then the leaf $\tilde{\cF}^{ss}_g(p)$ is $\tilde{g}$-invariant.
\end{claim}
\proof
Since $\cF^c$ is topologically Anosov and $\cF^{cs}$ is the stable foliation of the center flow, the foliation $\cF^c_g$ has the same feature.
The strong stable leaf $\tilde{\cF}^{ss}_g(p)$ separates the leaf  $\tilde{\cF}^{cs}(y)$, which is a plane, into two connected components $P_1$ and $P_2$ such that the center leaves converge  in $P_1$ and separate in $P_2$. Once again by the topological Anosov property of center foliation, one has that $\tilde{\cF}^{ss}_g(\tilde{\cF}^c_g(q))\subset P_1$.
By Claim~\ref{c.keep center leaf invariant},  every center leaf in $\tilde{\cF}^{cs}(y)$ is $\tilde{g}$-invariant, hence $\tilde{\cF}^{ss}_g(\tilde{\cF}^c_g(q))$ is a $\tilde{g}$-invariant set and the center leaf $\tilde{\cF}^c(p)$ is also $\tilde{g}$-invariant.
  Since $\tilde{g}$ sends a boundary leaf to a boundary leaf and  $\tilde{\cF}^{ss}_g(\tilde{\cF}^c_g(q))$ is a path connected invariant set,  one has that $\tilde{g}(\tilde{\cF}^{ss}(p))\subset P_1$.  The boundary leaf $\tilde{\cF}^{ss}_g(\tilde{g}(p))$ also separates the leaf  $\tilde{\cF}^{cs}(y)$ into two connected components $P^\prime_1$ and $P^\prime_2$ such that the center leaves converge  in $P^\prime_1$ and separate in $P^\prime_2$.
  Hence, one has that  $\tilde{\cF}^{ss}_g(\tilde{\cF}^c_g(q))\subset P^\prime_1$.
  By the invariance of $\tilde{\cF}^c(p)$, one has that $\tilde{\cF}^c(p)$ intersects $\tilde{\cF}^{ss}_g(\tilde{g}(p))$, and by transversality, the intersection is unique and we denote it as $z$. If $\tilde{\cF}^{ss}(p)$ is not $\tilde{g}$-invariant, then the connected component of $\tilde{\cF}^c(p)\backslash \{z\}$ which does not contain $p$ is contained in $P_2^{\prime}$. Also the connected component of  $\tilde{\cF}^c(p)\backslash \{p\}$ which does not contain $z$ is in $P_2$.  As a consequence, on the leaf $\tilde{\cF}^{cs}(y)$,   one has that the center leaf $\tilde{\cF}^c(p)$ is uniformly away from the center leaf $\tilde{\cF}^c(q)$, which contradicts to the  topologically Anosov property of the center foliation.
\endproof

 Since every center leaf in $\tilde{\cF}_g^{cs}(y)$ is  $\tilde{g}$-invariant and  every center leaf  intersects  $\tilde{\cF}_g^{ss}(p)$ in at most one point,
one has that every point in $\tilde{\cF}_g^{ss}(p)$ is a fixed point of $\tilde{g}$, contradicting to the fact that $\tilde{\cF}_g^{ss}(p)$ is a strong stable leaf,
ending the proof of Lemma~\ref{l.completeness in repeller}.
\endproof

Now, we consider the center stable leaves in the region $N\backslash h_g(\cR)$.
Assume, on the  contrary,  that  there exists a point $x$ such that
$\cF_g^{ss}(\cF_g^c(x))\subsetneq \cF_g^{cs}(x),$
 then let $p$ be a point  such that $\cF_g^{ss}(p)$ is a boundary leaf of $\cF_g^{ss}(\cF_g^c(x))$.
Since $h_g(\cA)$ is the maximal invariant  set in $N^-$ and is saturated by center unstable leaves,
 there exists  an integer $n$ large enough such that  $\cF_g^{ss}(g^n(p))$ intersects $h_g(\cA)$ in a point $q^\prime$.
Since   every center leaf in $h_g(\cA)$ is $g$-invariant, the center leaf through $q^{\prime}$ is contained in $\cF_g^{cs}(x)$ and intersects $\cF_g^{ss}(p)$. We denote by $q=g^{-n}(q^{\prime})$, then $q\in\cF_g^{ss}(p)\cap\cF_g^{c}(q^{\prime})$.

Now, we lift these leaves to the universal cover. Let $\tilde{x}$, $\tilde{p}$ and $\tilde{q}$ be the  lifts of $x$, $p$ and $q$ respectively
such that they are on the same center stable leaf $\tilde{\cF}_g^{cs}(\tilde{x})$ and $\tilde{\cF}^c_g(\tilde{q})$ intersects $\tilde{\cF}^{ss}_g(\tilde{p})$.
\begin{lemma}\label{l.fixed leaf} There exist a  lift $\tilde{g}$ of $g$ and  a
  center leaf $L\subset \tilde{\cF}_g^{cs}( \tilde{x} )$ such that
 \begin{itemize}\item the leaf $L$ is  disjoint from $\tilde{\cF}_g^{ss}(\tilde{p})$;
 \item the center leaves $\tilde{\cF}^c_g(\tilde{q})$ and $L$ are  $\tilde{g}$-invariant.
\end{itemize}
\end{lemma}
\proof  If $\cF_g^c(x)$ is $g$-invariant,  then at least  one of the  invariant center leaves $\cF^c_g(x)$ and $\cF^c_g(q)$ is not compact.
Assume that $\cF^c_g(q)$ is not compact (the other case follows analogously). Then there exists a unique lift $\tilde{g}$ of $g$ such that the center leaf $\tilde{\cF}^c_g(\tilde{q})$ is $\tilde{g}$-invariant.
By Theorem~\ref{thm.conjugate to Anosov flow} and  leaf conjugacy, one has the following:
\begin{itemize}
\item the strip bounded by $\tilde{\cF}^c_g(\tilde{q})$ and $\tilde{\cF}^c_g(\tilde{x})$ is trivially foliated by center leaves;
\item there exists a  $C^1$ curve  $\ell\subset \tilde{\cF}^c_g(\tilde{q})$ with infinity length such that for any point $z\in\ell$, the strong stable leaf
$\tilde{\cF}^{ss}_g(z) $ through $z$ intersects the center leaf  $\tilde{\cF}^c_g(\tilde{x})$.
\end{itemize}
Since  $f_b$ preserves the orientation of the center foliation, by leaf conjugacy, there exists a point $w\in\ell$ such that $\tilde{g}(w)\in\ell$.
 We take the strong stable segment $\tilde{\sigma}(t)$ through $w$ whose two endpoints are contained in  $\tilde{\cF}^c_g(\tilde{q})$ and $\tilde{\cF}^c_g(\tilde{x})$ respectively. Now, one can check that the arguments in Claim~\ref{c.keep center leaf invariant} can be applied and one gets that $\tilde{\cF}_g^c(\tilde{x})$ is $\tilde{g}$-invariant.  We only need to  take
$L=\tilde{\cF}_g^c(\tilde{x})$.

  If center leaf $\cF_g^c(x)$ is not $g$-invariant, then the center leaf $h^{-1}_g(\cF_g^c(x))$ for $f_b$ intersects  the transverse torus  $T_1$  and it is not  $f_b$-invariant.  Consider the  the connected component $\cP$ of  $\cF^{cs}(h^{-1}_g(x))\backslash \cA$ which contains the center leaf $\cF^c(h^{-1}_g(x))$, then $\cP$ is a topological plane. Recall that on each transverse torus, the foliation $\cF^{cs}$ induces a foliation consisting of  exactly  two Reeb components.
  Since $\cF^c(h^{-1}_g(x))$ is not $f$-invariant, by Lemma~\ref{l.existence of fixed center leaf},  the intersection between $\cP$ and $T_1$ can't be a circle, hence the intersection is  a line $\ell_1$;  we identify $\ell_1$ with $\mathbb{R}$, then $\ell_1$ accumulates to two circles $S_1,S_2$ when it tends to infinity; moreover, the circles $S_1,S_2$ are contained in the stable manifolds of two different periodic orbits of  $\psi_t$.
  Hence   the boundary of $\cP$ restricted to $\cF^{cs}(h^{-1}_g(x))$  consists of  two center leaves belonging to the unstable manifolds of  two different periodic orbits of $\psi_t$ in $\cA$.
  Then, by leaf conjugacy,  there exist two  leaves $L_1, L_2$ of the foliation $\tilde{\cF}_g^c$ such that
  \begin{itemize}
  \item $L_1\cup L_2\subset \tilde{\cF}_g^{cs}(\tilde{x})$;
   \item the leaves $L_1$ and $L_2$ bound a strip  containing $\tilde{\cF}^c(\tilde{x})$;
   \item  the projections of   the center leaves $L_1, L_2$ on the base manifold are $g$-invariant.
    \end{itemize}
   By the choices of $L_1,L_2$, neither $\pi(L_1)$ nor $\pi(L_2)$ is a compact leaf,
   hence one can apply the argument in the first case and one gets a lift $\tilde{g}$ of $g$ such that the center leaves $L_1,L_2, \tilde{\cF}_g^c(\tilde{q})$ are $\tilde{g}$-invariant.
 Then there exists $i\in\{1,2\}$ such that $\tilde{\cF}_g^c(\tilde{x})$ is contained in the strip $S$ bounded by $L_i$ and $\tilde{\cF}_g^c(\tilde{q})$.
Since $\tilde{\cF}_g^{ss}(\tilde{p})$ does not intersect $\tilde{\cF}_g^c(\tilde{x})$, the strong stable leaf
$\tilde{\cF}_g^{ss}(\tilde{p})$ does not intersect   $L_i$.
We take $L=L_i$, ending the proof of Lemma~\ref{l.fixed leaf}.
\endproof

\begin{claim}\label{c.L}There exists a $\tilde{g}$-fixed point on $\tilde{\cF}_g^c(\tilde{q})$ whose strong stable leaf does not intersect $L$.
\end{claim}
\proof If $\tilde{q}$ is a $\tilde{g}$-fixed point, we are done. Now, we assume that  $\tilde{q}$ is not a fixed point.
Recall that $\tilde{q}\in\widetilde{h_g(\cA)}$.
 We denote by  $I_q$  the connected component of $\tilde{\cF}_g^{c}(\tilde{q})\backslash \{\tilde{q}\}$ such that under leaf conjugacy,
 it corresponds to the forward orbit of $\tilde{h}_g^{-1}(q)$ under the Anosov flow $\psi_t$.  Up to replacing $\tilde{g}$ by $\tilde{g}^{-1}$, we can assume that $\tilde{g}(\tilde{q})$ is contained in the interior of $I_q$. Since  the leaf $L$ is fixed by $\tilde{g}$, the strong stable leaf through the orbit of $\tilde{q}$ is disjoint from $L$. We identify $I_q$ with $(0,+\infty)$, where $q$ corresponds to $0$.
 By Theorem~\ref{thm.conjugate to Anosov flow}, for the points on $I_q$ tending to infinity, their  strong stable leaves would intersect $L$. Hence, one has that the forward orbit of $\tilde{q}$ tends to a $\tilde{g}$-fixed point whose strong stable leaf is disjoint from $L$.
\endproof

By  Claim~\ref{c.L}, for notational convenience, one can assume that   $\tilde{q}$ is the fixed point of $\tilde{g}$ in $\tilde{\cF}^c(\tilde{q})$. Since $\cF_g^c(q)$  is contained in $h_g(\cA)$  and  $\cA$ is an attractor of the Anosov flow $\psi_t$,
restricted to the metric on the center stable leaf $\cF_g^{cs}(q)$,
the leaf   $\cF_g^c(q)$ is accumulated by $g$-invariant center leaves.
Hence, there exist center leaves on $\tilde{\cF}_g^{cs}(\tilde{q})$ contained in $\widetilde{h_g(\cA)}$
which intersect the strong stable leaf $\tilde{\cF}_g^{ss}(\tilde{q})$. Following the argument in Claim~\ref{c.keep center leaf invariant},
one can check that these center leaves are invariant under $\tilde{g}$. Recall that in $\tilde{\cF}^{cs}(\tilde{q})$, each center leaf intersects each strong stable leaf in at most one point.
Since  these center leaves and the strong stable leaf $\tilde{\cF}_g^{ss}(\tilde{q})$ are $\tilde{g}$-invariant,
one has that the intersection between these center leaves and  $\tilde{\cF}_g^{ss}(\tilde{q})$ consists of $\tilde{g}$-fixed points which contradicts to the uniform contraction of $\tilde{g}$ along $\tilde{\cF}_g^{ss}(\tilde{q})$. Hence, the center stable foliation of $g$ is complete.

Similar argument applies for the center unstable foliation of $g$.
 \endproof

We point out that we use ~\cite{BPP} example to get robust completeness instead of the general ~\cite{BZ} example, because we need the orbits (for Anosov flow) in the attracting set (resp. repelling set) are not isolated on it stable (resp. unstable) manifolds.

\vspace{5mm}
{\bf Acknowledgement:}
I would like to thank Christian Bonatti whose questions motivated this paper. I  thank Rafael Potrie  for many helpful suggestions on writing and useful comments on math,  and Shaobo Gan for listening the proof of this paper.  I also thank Lan Wen, Sylvain Crovisier and Dawei Yang for useful comments. Finally, I want to thank ERC Advanced Grant NUHGD for financial support.

\bibliographystyle{plain}


\vspace{2mm}

\noindent Jinhua Zhang,

\noindent{\small Laboratoire de Math\'ematiques d'Orsay,\\
UMR 8628 du CNRS,\\
Universit\'e Paris-Sud,
 91405 Orsay, France.}
\\
\noindent {\footnotesize{E-mail :  zjh200889@gmail.com $\&$ jinhua.zhang@u-psud.fr}}\\

\end{document}